\documentclass[aap,MSNbibl,nameyear,seceqn,dvips]{arximspdf}
\usepackage{graphics}
\doi{10.1214/09-AAP631}
\volume{20}
\issue{2}
\pubyear{2010}
\firstpage{753}
\lastpage{783}

\makeatletter
\def\eqref#1{(\ref{#1})}
\DeclareMathAlphabet\mathcaligr{OMS}{cmsy}{m}{n}
\newtheorem{lemma}{Lemma}
\newtheorem{theorem}{Theorem}

\newproclaim{category}{Category}
\newproclaim{definition}{Definition}
\makeatother

\begin{document}
\begin{frontmatter}

\title{Positive recurrence of reflecting Brownian motion in three
dimensions}
\runtitle{Reflecting Brownian motion}

\begin{aug}
\author[A]{\fnms{Maury} \snm{Bramson}\ead[label=e1]{bramson@math.umn.edu}},\thanksref{a1}
\author[B]{\fnms{J. G.} \snm{Dai}\thanksref{a2}\ead[label=e2]{dai@gatech.edu}} \and
\author[C]{\fnms{J. M.} \snm{Harrison}\ead[label=e3]{harrison\_michael@gsb.stanford.edu}\corref{}}
\runauthor{M. Bramson, J. G. Dai and J. M. Harrison}
\affiliation{University of Minnesota, Georgia Institute of Tecnology
and Stanford University}
\address[A]{M. Bramson\\ School of Mathematics\\ University of
Minnesota\\ Minneapolis, Minnesota 55455\\ USA\\\printead{e1}} 
\address[B]{J. G. Dai\\ H. Milton Stewart School of Industrial\\ \quad
and Systems Engineering\\ Georgia Institute of Technology\\
Atlanta, Georgia 30332\\ USA\\\printead{e2}}
\address[C]{J. M. Harrison\\ Graduate School of Business\\ Stanford
University\\ Stanford, California 94305-5015\\USA\\\printead{e3}}
\thankstext{a1}{Supported in part by NSF Grant CCF-0729537.}
\thankstext{a2}{Supported in part by NSF Grants CMMI-0727400 and
CMMI-0825840, and by an IBM Faculty Award.}
\end{aug}

\received{\smonth{11} \syear{2008}}
\revised{\smonth{7} \syear{2009}}

%
\begin{abstract}
Consider a semimartingale reflecting Brownian motion (SRBM) $Z$ whose state
space is the $d$-dimensional nonnegative orthant. The~data for
such a process are a drift vector $\theta$, a nonsingular $d\times d$
covariance matrix $\Sigma$, and a $d\times d$ reflection matrix $R$ that
specifies the boundary behavior of $Z$. We say that $Z$ is positive
recurrent, or stable, if the expected time to hit an arbitrary open
neighborhood of the origin is finite for every starting state.

In dimension $d = 2$, necessary and sufficient conditions for
stability are known, but fundamentally new phenomena arise in higher
dimensions. Building on prior work by {El Kharroubi}, {Ben Tahar} and
Yaacoubi [\textit{Stochastics Stochastics Rep.} \textbf{68} (2000) 229--253,
\textit{Math. Methods Oper. Res.} \textbf{56} (2002) 243--258], we
provide necessary and
sufficient conditions for stability of SRBMs in three dimensions; to
verify or refute these conditions is a simple computational task.
As a byproduct, we find that the fluid-based criterion of
Dupuis and Williams [\textit{Ann. Probab.} \textbf{22} (1994)
680--702] is not only sufficient but also necessary for
stability of SRBMs in three dimensions. That is, an SRBM in three
dimensions is positive recurrent if and only if every path of the
associated fluid model is attracted to the origin. The~problem of
recurrence classification for SRBMs in four and higher dimensions
remains open.
\end{abstract}

%

\end{frontmatter}
%

\section{Introduction}\label{sec1}

This paper is concerned with the class of $d$-dimensional
diffusion processes called \textit{semimartingale reflecting
Brownian motions} (SRBMs), which arise as approximations for open
$d$-station queueing networks of various kinds; cf. \citet{harngu93}
and \citeauthor{wil95} (\citeyear{wil95}, \citeyear{wil96}). The state space for a
process $Z =\{ Z(t),   t \ge0 \}$ in this class is $S=\mathbb{R}^d_+$
(the nonnegative orthant). The~data of the process are a drift
vector $\theta$, a nonsingular covariance matrix $\Gamma$, and a
$d\times d$ ``reflection matrix'' $R$ that specifies boundary
behavior. In the interior of the orthant, $Z$ behaves as an
ordinary Brownian motion with parameters $\theta$ and $\Gamma$, and
roughly speaking, $Z$ is pushed in direction $R^{j}$ whenever the
boundary surface {$\{z\in S\dvtx  z_j=0\}$} is hit, where $R^{j}$ is
the $j$th column of $R$, for $j = 1, \ldots, d$.
To make this description
more precise, one represents $Z$ in the form
%
\begin{eqnarray}
\label{11}
&& Z(t)=X(t) + RY(t), \qquad  t \ge0,
\end{eqnarray}
where $X$ is an unconstrained Brownian motion with drift vector $\theta$,
covariance matrix $\Gamma$, and $Z(0) = X(0)\in S$, and $Y$ is
a $d$-dimensional process with components $Y_{1}$ , \ldots, $Y_{d}$
such that
%
\begin{eqnarray}
\label{12}
&&\qquad  Y \mbox{ is continuous and nondecreasing with } Y(0) = 0, \\
\label{13}
&&\qquad  Y_{j} \mbox{ only increases at times $t$ for which }  Z_{j}(t)
= 0, \quad   j = 1,
\ldots, d,\quad  \mbox{and} \\
&&\qquad  \label{14} Z(t)\in S, \qquad   t \ge0.
\end{eqnarray}
The~complete definition and essential properties of the diffusion
process $Z$ will be reviewed in Appendix~\hyperref[sec:srbm]{A}, where we
also discuss the notion of positive recurrence. As usual in Markov
process theory, the complete definition involves a family of
probability measures $\{\mathbb{P}_{x},  x\in S\}$ that specify the
distribution of $Z$ for different starting states; informally, one
can think of $\mathbb{P}_{x}(\cdot)$ as a conditional probability
given that $Z(0) = x$. Denoting by $\mathbb{E}_{x}$ the expectation
operator associated with $\mathbb{P}_{x}$ and setting $\tau
_{A}=\inf\{t\ge0\dvtx  Z(t)\in A \}$, we say that $Z$ is positive
recurrent if
\mbox{$\mathbb{E}_{x}(\tau_{A}) < \infty$} for any $x\in S$ and any
open neighborhood $A$ of the origin (see Appendix~\hyperref[sec:srbm]{A} for
elaboration). For ease of expression, we use the terms ``stable''
and ``stability'' as synonyms for ``positive recurrent'' and
``positive recurrence,'' respectively.

In the foundational theory for SRBMs, the following classes of
matrices are of interest. First, a $d\times d$ matrix $R$ is said to
be an $\mathcaligr{S}$-\textit{matrix} if there exists a $d$-vector $w
\ge$ 0 such that $Rw>0$ (or equivalently, if there exists $w > 0$
such that $Rw > 0$), and $R$ is said to be
\textit{completely}-$\mathcaligr{S}$ if each of its principal
submatrices is an $\mathcaligr{S}$-matrix. (For a vector $v$, we write
$v>0$ to mean that each component of $v$ is positive, and we write
$v\ge0$ to mean that each component of $v$ is nonnegative.) Second,
a square matrix is said to be a $\mathcaligr{P}$-matrix if all of its
principal minors are positive (that is, each principal submatrix of
$R$ has a positive determinant). \mbox{$\mathcaligr{P}$-matrices} are a
subclass of completely-$\mathcaligr{S}$ matrices; the still more
restrictive class of $\mathcaligr{M}$-matrices is defined as in Chapter 6
of \citet{berple79}. References for the following key results can be
found in the survey paper by \citet{wil95}: there exists a diffusion
process $Z$ of the form described above if and only if $R$ is a
completely $\mathcaligr{S}$ matrix; and moreover, $Z$ is unique in
distribution whenever it exists.

Hereafter we assume that $R$ is completely-$\mathcaligr{S}$. Its diagonal
elements must then
be strictly positive, so we can (and do) assume without loss of generality
that
%
\begin{equation}
\label{eq:new1.5}
R_{ii} = 1\qquad  \mbox{for all } i = 1, \ldots, d.
\end{equation}
This convention is standard in the SRBM literature; in Sections
\ref{sec:spiraling} through \ref{sec:lemmas} of this paper (where our
main results are proved) another convenient normalization of problem
data will be used. Appendix \hyperref[sec:normalization]{B} explains the
scaling procedures that justify both (\ref{eq:new1.5}) and the
normalized problem format assumed in Sections \ref{sec:spiraling}
through \ref{sec:lemmas}.

We are concerned in this paper with conditions that assure the
stability of
$Z$. An important condition in that regard is the following:
%
\begin{equation}\label{15}
R \mbox{ is nonsingular}\quad  \mbox{and}\quad R^{-1}\theta< 0.
\end{equation}
If $R$ is
an $\mathcaligr{M}$-matrix, then (\ref{15}) is known to be necessary
and sufficient for stability of $Z$; \citet{harwil87a} prove that
result and explain how the \mbox{$\mathcaligr{M}$-matrix} structure arises
naturally in queueing network applications.

\citet{ElKharroubiEtAl00} further prove the following three results:
first, (\ref{15}) is necessary for stability in
general; second, when $d=2$, one has stability if and only if
(\ref{15}) holds and $R$ is a $\mathcaligr{P}$-matrix; and third, (\ref
{15}) is not sufficient for stability in three
and higher dimensions, even if $R$ is a $\mathcaligr{P}$-matrix. In
Appendix \hyperref[sec:proofNecessa]{C} of this paper, we provide an
alternative proof for the first of these results, one that is much
simpler than the original proof by \citet{ElKharroubiEtAl00}.
Appendix A of \citet{HarrisonHasenbein08} contains an alternative proof
of the second result. Section \ref{sec:bek-example} of this paper
reviews the ingenious example by \citet{BernardElKharroubi91} that
serves to establish the third result, that
an SRBM can be unstable, cycling to infinity even if
(\ref{15}) holds; Theorem~\ref{thm:4} in this paper, together with
the examples provided in Section \ref{sec:divergent}, shows
that instability can also occur in other ways when (\ref{15}) holds.

A later paper by \citet{ElKharroubiEtAl02} established sufficient
conditions for stability of SRBMs in three dimensions, relying heavily
on the foundational theory developed by \citet{dupwil94}.
In this paper, we show that the conditions identified by\break
\citet{ElKharroubiEtAl02} are also necessary for stability when $d=3$; the
relevant conditions are easy to verify or refute via simple
computations. As a complement to this work, an alternative proof of
the sufficiency result by \citet{ElKharroubiEtAl02} is also being
prepared for submission; cf.~\citet{DaiHarrison09b}.

The~remainder of the paper is structured as follows. First, to allow
precise statements of the main results, we introduce in Section
\ref{sec:lcp} the ``fluid paths'' associated with an SRBM, and the
linear complementarity problem that arises in conjunction with linear
fluid paths. That section, like the paper's first three appendices,
considers a general dimension $d$, whereas all other sections in the
body of the paper consider $d=3$ specifically. Section
\ref{sec:bek-example} identifies conditions under which fluid paths
spiral on the boundary of the state space $S$. Section~\ref
{sec:3d-results} states our main conclusions, which are achieved
by combining the positive results of \citet{ElKharroubiEtAl02} with
negative results that are new; Figure~\ref{fig:3} in Section~\ref
{sec:3d-results} summarizes succinctly the necessary and
sufficient conditions for stability when $d=3$, and indicates
which components of the overall argument are old and which are new.
In Sections~\ref{sec:spiraling} through \ref{sec:lemmas}, we prove the
new ``negative results'' referred to above, dealing first with the
case where fluid paths spiral on the boundary, and then with the case
where they do not. As stated above, Appendix~\hyperref[sec:srbm]{A} reviews
the precise definition of SRBM; Appendix~\hyperref[sec:normalization]{B}
explains the scaling procedures that give rise to normalized problem
formats and Appendix~\hyperref[sec:proofNecessa]{C} contains a relatively
simple proof that (\ref{15}) is necessary for stability. Finally,
Appendix~\hyperref[sec:AppendixD]{D} contains several technical lemmas that are
used in the probabilistic arguments of Section~\ref{sec:lemmas}.

\section{Fluid paths and the linear complementarity problem}
\label{sec:lcp}

\begin{definition}
A fluid path associated with the data $(\theta, R)$ is a pair of
continuous functions $y, z\dvtx  [0,\infty)\to\mathbb{R}^d$ that satisfy
the following
conditions:
%
\begin{eqnarray}
\label{eq:fl-1}
&& z(t) = z(0) + \theta t + Ry(t)\qquad  \mbox{for all } t \ge0, \\
&& z(t) \in S\qquad  \mbox{for all } t \ge0,\label{eq:fl-2} \\
&& y(\cdot) \mbox{ is continuous and nondecreasing with } y(0) =
0,\label{eq:fl-3} \\
\label{eq:fl-4}&& y_j(\cdot) \mbox{ only increases when } z_j(\cdot) = 0, \mbox{ i.e., }\nonumber\\[-8pt]\\[-8pt]
&&\int_0^\infty z_j(t) \, d y_j(t)=0, \qquad  (j =
1,\ldots, d).\nonumber
\end{eqnarray}
\end{definition}

\begin{definition}
We say that a fluid path $(y, z)$ is \textit{attracted to the origin}
if $z(t)
\to0$ as $t \to\infty$.
\end{definition}

\begin{definition}
A fluid path $(y, z)$ is said to be \textit{divergent} if $|z(t)| \to
\infty$
as $ t \to\infty$, where, for a vector $u=(u_i)\in
\mathbb{R}^d$, $|u|=\sum_{i} |u_i|$.
\end{definition}

\begin{theorem}[\mbox{[}Dupuis and Williams (\citeyear{dupwil94})\mbox{]}] Let $Z$ be a $d$-dimensional
SRBM with data
$(\theta,\Gamma, R)$. If every fluid path associated with $(\theta,
R)$ is
attracted to the origin, then $Z$ is positive recurrent.
\label{thm:2}
\end{theorem}

\begin{definition}
A fluid path $(y, z)$ is said to be \textit{linear} if it has the form $y(t)
= u t$ and $z(t) = v t$, $t \ge0$, where $u,v
\ge0$.
\end{definition}

 Linear fluid paths are in one-to-one correspondence
with solutions of the following linear complementarity problem
(LCP): Find vectors $u=(u_i)$ and $v=(v_i)$ in $\mathbb{R}^d$ such that
%
\begin{eqnarray}
 u, v&\ge&0, \label{eq:LCP1} \\
 v&=&\theta+Ru, \label{eq:LCP2} \\
 u\cdot v& =& 0, \label{eq:LCP3}
\end{eqnarray}
where $u\cdot v=\sum_{i}u_i v_i$ is the inner product of $u$ and
$v$.
[See \citet{cps92} for a systematic account of the theory associated
with the general problem (\ref{eq:LCP1})--(\ref{eq:LCP3}).]

\begin{definition}
A solution $(u,v)$ of the LCP is said to be \textit{stable} if
$v=0$ and to be \textit{divergent} otherwise. It is said to be
\textit{nondegenerate} if $u$ and $v$ together have
exactly $d$ positive
components, and to be \textit{degenerate} otherwise. A stable, nondegenerate
solution of the LCP is called \textit{proper}.
\end{definition}

\begin{lemma}\label{lem:lem1}
Suppose that (\ref{15}) holds. Then $(u^*,0)$ is a proper solution of
the LCP, where
%
\begin{equation}
\label{eq:ustar}
u^* = -R^{-1}\theta,
\end{equation}
and any other solution of the LCP must be divergent.
\end{lemma}

\begin{pf}
The~first statement is obvious. On the other hand, for any stable
solution $(u,0)$ of
the LCP, we have from (\ref{eq:LCP2}) that $\theta+Ru = 0$; since
(\ref{15}) includes the requirement that $R$ be nonsingular, $u
=- R^{-1}\theta= u^*$. That is, there cannot exist a stable solution other
than $(u^*,0)$, which is equivalent to the second statement of the
lemma.
\end{pf}

\begin{figure}[b]

\includegraphics{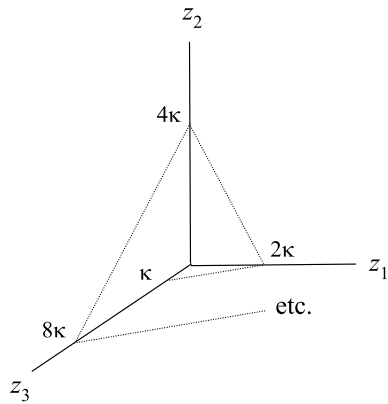}

\caption{Fluid model behavior of the B{\&}EK example.}
\label{fig2}
\end{figure}

\section{Fluid paths that spiral on the boundary}\label{sec:bek-example}

\citet{BernardElKharroubi91} devised the following ingenious example
with $d=3$, referred to hereafter as the B{\&}EK example:
let
\[
\theta=
\pmatrix{
{-1}  \cr
{-1}  \cr
{-1}  }\quad   \mbox{and}\quad    R=\pmatrix{
1  & 3  & 0  \cr
0  & 1  & 3  \cr
3  & 0  & 1  }.
\]
This
reflection matrix $R$ is completely $\mathcaligr{S}$ (moreover, it is
a $\mathcaligr{P}$-matrix), so $Z$ is a well-defined SRBM.
(The~covariance matrix $\Gamma$ is immaterial to the\break discussion that
follows, provided only that it is nonsingular.) As\break
\citet{BernardElKharroubi91} observed, the unique fluid path with these
process data, starting from $z(0) = (0,0,\kappa)$ with $\kappa> 0$,
is the one pictured in Figure \ref{fig2}; it travels in a
counter-clockwise and piecewise linear fashion on the boundary, with
the first linear \mbox{segment} ending at $(2\kappa,0,0)$, the second one
ending at $(0,4\kappa,0)$, and so forth. \citet{ElKharroubiEtAl00}
proved that an SRBM with these data is not stable, showing that if
$\kappa$ is large then $|Z(t)-z(t)|$ remains forever small (in a
certain sense) with high probability.

To generalize the B{\&}EK example, let $C_{1}$ be the set of ($\theta,R)$
pairs that satisfy the following system of inequalities [here $R_{ij}$
denotes the ($i,j)${th} element of $R$, or equivalently, the $i${th}
element of
the column vector $R^{j}$]:
%
\begin{eqnarray}
\label{eq:31}
 \theta&<&0, \\
 \theta_{1} &>& \theta_{2}R_{12}\quad   \mbox{and}\quad   \theta
_{3} <\theta
_{2}R_{32}, \label{eq:32}\\
 \theta_{2} &>& \theta_{3}R_{23}\quad  \mbox{and}\quad   \theta_{1}
< \theta
_{3}R_{13 },\label{eq:33} \\
\theta_{3} &>& \theta_{1}R_{31} \quad  \mbox{and}\quad   \theta
_{2}< \theta
_{1}R_{21}.
\label{eq:34}
\end{eqnarray}
[Notation used in this section agrees with that
of El~Kharroubi, Tahar and
Yaacoubi (\citeyear{ElKharroubiEtAl00}, \citeyear{ElKharroubiEtAl02}) in all essential respects,
but is different in a few minor
respects.]

To explain the meaning of these inequalities, we consider a fluid
path associated with $(\theta,R)$ that starts from $z(0) = (0,0,\kappa)$,
where $\kappa>0$; it is the unique fluid path starting from that state,
but that fact will not be used in our formal results. Over an initial time
interval $[0,\tau_{1}]$, the fluid path is linear and adheres to the
boundary $\{z_{2} = 0\}$, as in Figure \ref{fig2}. During that
interval one has
$\dot{y}(t)=(0,- \theta_2 ,0{)}'$
and hence the fluid path has the constant velocity vector
%
\begin{equation}
\label{eq:velocity}
\dot z(t) = \theta+ R \dot y(t) = \theta
-\theta_2 R^2 =
\pmatrix{
\theta_1 -\theta_2 R_{12} \vspace*{2pt}\cr
0 \vspace*{2pt}\cr
\theta_3 -\theta_2 R_{32}}
.
\end{equation}
%

Thus (\ref{eq:31}) and (\ref{eq:32}) together give the following: as
in Figure \ref{fig2}, a fluid
path starting from state $(0,0,\kappa)$ has an initial linear segment in
which $z_{3}$ decreases, $z_{1 }$~increases and $z_{2 }$ remains at
zero; that
initial linear segment terminates at the point $z(\tau_{1})$ on the
$z_{1 }$ axis that has
\[
z_1 ( \tau_1 )=\biggl( {\frac{ \theta_1 - \theta_2 R_{12} }{ \theta
_2 R_{32} - \theta_3 }} \biggr)\kappa> 0.
\]

Similarly, from (\ref{eq:31}) and (\ref{eq:33}), the fluid path is
linear over an ensuing time interval $[\tau_{1},\tau_{2}]$, with
$z_{1}$ decreasing, $z_{2 }$ increasing and $z_{3 }$ remaining at
zero; that second linear segment terminates at the point $z(\tau
_{2})$ on the $z_{2 }$ axis that has
\[
z_2 ( \tau_2 )=\biggl( {\frac{ \theta_1 - \theta_2 R_{12} }{ \theta
_2 R_{32} - \theta_3 }} \biggr)\biggl( {\frac{ \theta_2 -
\theta_3 R_{23} }{ \theta_3 R_{13} - \theta_1 }}
\biggr)\kappa> 0.
\]
Finally, from (\ref{eq:31}) and (\ref{eq:34}), the fluid path is
linear over a next time
interval $[\tau_{2},\tau_{3}]$, with $z_{2}$ decreasing, $z_{3
}$ increasing and $z_{1 }$ remaining at zero; that third linear segment
terminates at the point $z(\tau_{3})$ on the $z_{3 }$ axis that has
$z_{3}(\tau_{3})= \beta_1 (\theta,R)\kappa$,
where
%
\begin{equation}
\label{eq:35}
\beta_1 (\theta,R)=\biggl( {\frac{ \theta
_1 - \theta_2 R_{12} }{
\theta_2 R_{32} - \theta_3 }}
\biggr)\biggl( {\frac{ \theta_2 - \theta_3
R_{23} }{ \theta_3 R_{13}
- \theta_1 }} \biggr)\biggl( {\frac{ \theta
_3 - \theta_1 R_{31} }{
\theta_1 R_{21} - \theta_2 }}
\biggr)> 0.
\end{equation}
Thereafter, the piecewise linear fluid path continues its
counter-clockwise spiral on the boundary in a self-similar fashion,
like the path pictured in Figure \ref{fig2}, except that in the
general case defined by (\ref{eq:31}) through (\ref{eq:34}), the spiral
may be either inward or outward, depending on whether $ \beta_1
(\theta,R)< 1$ or $ \beta_1 (\theta, R)> 1$.

To repeat, $C_{1}$ consists of all $(\theta,R)$ pairs that satisfy
(\ref{eq:31}) through (\ref{eq:34}), and the \textit{single-cycle
gain} $ \beta_1 (\theta,R)$ for such a pair is defined by
(\ref{eq:35}). As we have seen, fluid paths associated with problem
data in $C_{1}$ spiral counter-clockwise on the boundary of $S$. Now
let $C_{2}$ consist of all $(\theta,R)$ pairs that satisfy (\ref{eq:31})
and further satisfy (\ref{eq:32}) through (\ref{eq:34}) \textit{with
all six of the
strict inequalities reversed}. It is more or less obvious that
$(\theta,R)$ pairs in $C_{2}$ are those giving rise to
\textit{clockwise} spirals on the boundary, and the appropriate analog
of (\ref{eq:35}) is
\begin{eqnarray}
\label{eq:36}
\beta_2 (\theta,R)&=&\frac{1}{ \beta
_1 (\theta,R)}\nonumber\\[-8pt]\\[-8pt]
&=&\biggl( {\frac{ \theta_3 -
\theta_2 R_{32} }{ \theta_2
R_{12} - \theta_1 }} \biggr)\biggl(
{\frac{ \theta_1 - \theta_3
R_{13} }{ \theta_3 R_{23}
- \theta_2 }} \biggr)\biggl( {\frac{ \theta
_2 - \theta_1 R_{21} }{
\theta_1 R_{31} - \theta_3 }}
\biggr)> 0.\nonumber
\end{eqnarray}
Hereafter we define $C= C_1 \cup C_2 $, $\beta(\theta,R)= \beta_1
(\theta,R)$ for $(\theta,R)\in C_1 $ and $\beta(\theta,R)= \beta
_2 (\theta,R)$ for $(\theta,R)\in C_2 $. Thus $C$ consists of all
$(\theta,R)$ pairs whose associated fluid paths spiral on the
boundary, and $\beta(\theta,R)$ is the single-cycle gain for such a
pair.

\section{Summary of results in three dimensions}
\label{sec:3d-results}

Theorem~\ref{thm:ELThm1} below is a slightly weakened version of
Theorem 1 by
\citet{ElKharroubiEtAl02}, which the original authors express in a more
elaborate notation; we have deleted one part of their result that is
irrelevant for current purposes. The~corollary that follows is
immediate from
Theorem 1 above (the Dupuis--Williams fluid stability criterion) and
Theorem \ref{thm:ELThm1}.

\begin{theorem}[\mbox{[}\citet{ElKharroubiEtAl02}\mbox{]}] \label{thm:ELThm1} Suppose
that (\ref{15}) holds and that either of the following additional
hypotheses is satisfied: \textup{(a)}~$(\theta,R)\in C$ and $\beta(\theta
,R) < 1$; or \textup{(b)} $(\theta,R)\notin C$ and the linear
complementarity problem (\ref{eq:LCP1})--(\ref{eq:LCP3}) has a unique
solution, which is the proper solution $(u^*,0)$ defined in
(\ref{eq:ustar}). Then all fluid paths associated $(\theta,R)$ are
attracted to the origin.
\end{theorem}

\begin{corollary*}
Suppose that (\ref{15}) holds and, in addition, either \textup{(a)} or \textup{(b)}
holds. Then $Z$ is positive recurrent.
\end{corollary*}

The~proof of Theorem~2 in \citet{ElKharroubiEtAl02} is not entirely
rigorous, containing verbal passages that mask significant technical
difficulties; an alternative proof that uses a linear Lyapunov
function to prove stability is given in \citet{DaiHarrison09b}.

\begin{figure}

\includegraphics{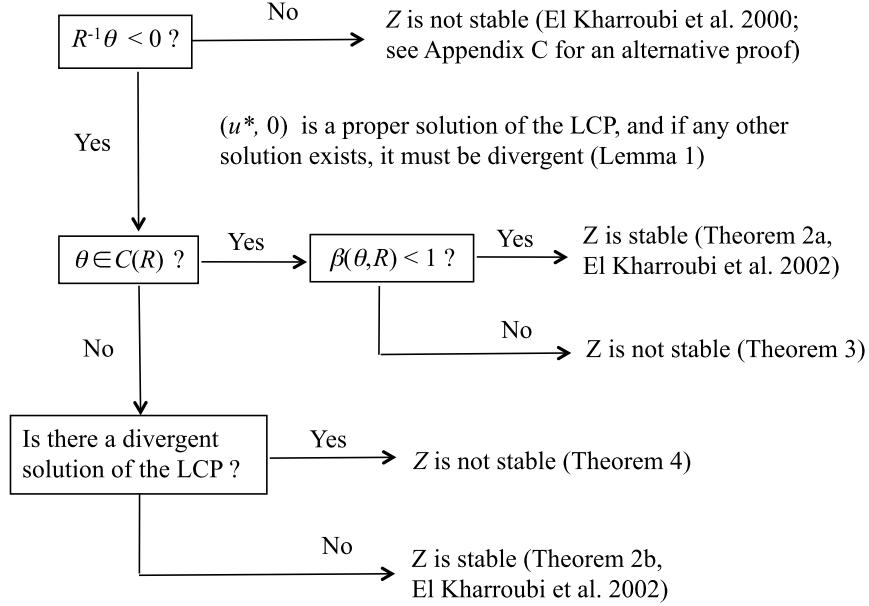}

\caption{Summary of results in three dimensions.}
\label{fig:3}
\end{figure}

The~new results of this paper are Theorems \ref{thm:3} and \ref{thm:4}
below, which will be proved in Sections \ref{sec:spiraling} through
\ref{sec:lemmas}.
Figure \ref{fig:3} summarizes the
logic by which these new results combine with previously known results
to provide necessary and sufficient conditions for stability (i.e.,
positive recurrence) of $Z$.

\begin{theorem}\label{thm:3}
If $(\theta,R)\in C$ and $\beta(\theta,R) \ge1$, then $Z$ is not
positive recurrent.
\end{theorem}

\begin{theorem}\label{thm:4}
Suppose that (\ref{15}) is satisfied. If there
exists a divergent solution for the linear complementarity problem
(\ref{eq:LCP1})--(\ref{eq:LCP3}), then $Z$ is not positive recurrent.
\end{theorem}

\section{Proof of Theorem~\protect\ref{thm:3}}
\label{sec:spiraling}

Throughout this section and the next, we assume without loss of
generality that our problem data satisfy not only (\ref{eq:new1.5})
but also
%
\begin{equation}
\label{eq:51}
\theta_i \in\{-1, 0, 1\}\qquad   \mbox{for } i = 1, 2, 3.
\end{equation}
Appendix \hyperref[sec:normalization]{B} explains the scaling procedures that
yield this normalized form. To prove Theorem~\ref{thm:3} we will
assume that $(\theta,R)\in C_1$ and $\beta_1(\theta,R)\ge1$, then
show that $Z$ is not stable; the proof of instability when $(\theta
,R) \in C_2$ and $\beta_2(\theta, R )\ge1$ is identical. Given the
normalizations (\ref{eq:new1.5}) and (\ref{eq:51}), the conditions
(\ref{eq:31}) through (\ref{eq:34}) that define $C_1$ can be restated
as follows:
%
\begin{eqnarray}
\label{eq:52}
 \theta&=&
(-1, -1, -1)' ,\\
 R_{12}, R_{23}, R_{31} &>& 1\quad   \mbox{and}\quad   R_{13}, R_{21},
R_{32} <1 .
\label{eq:53}
\end{eqnarray}

Let us now define a $3\times3$ matrix $V$ by setting $V_{ij} = R_{ij}
-1$ for $i, j = 1, 2, 3$. Then $V^j$ (the $j$th column of $V$) is the
vector $\theta- \theta_jR^j$ for $j=1, 2, 3$. Note that $V^2$ was
identified in
(\ref{eq:velocity}) as the velocity vector on the face $\{Z_2 = 0\}$
for a fluid path corresponding to $(\theta,R)$.

\begin{lemma}\label{lem:51}
Under the assumption that $\beta_1(\theta,R )\ge1$, there exists a
vector \mbox{$u > 0$} such that $u' V \ge0$, or equivalently, $u' V^j\ge
0$ for each $j = 1, 2, 3$.
\end{lemma}

\begin{pf}
From (\ref{eq:new1.5}) and (\ref{eq:53}) we have that
%
\begin{equation}
\label{eq:54}
V=
\pmatrix{  0 & a_2 & -b_3 \cr
-b_1 & 0 & a_3 \cr
a_1 & -b_2 & 0
}
,
\end{equation}
where $a_i, b_i > 0$ for $i = 1, 2, 3$. In this notation, the
definition (\ref{eq:35}) is as follows:
%
\begin{equation}
\label{eq:55}
\beta_1(\theta, R) = \frac{a_1a_2a_3}{b_1b_2b_3}.
\end{equation}
Setting
\[
u_1 = 1,\qquad   u_2 = \frac{a_1a_2}{b_1 b_2} \quad  \mbox{and}\quad
u_3 = \frac{a_2}{b_2},
\]
it is easy to verify that $u' V^1 = u' V^2 = 0$, and $u' V^3 = b_3(
\frac{a_1a_2a_3}{b_1b_2b_3} - 1)$. The~definition~\eqref{eq:55} and
our assumption that $\beta_1(\theta,R)\ge1$ then give $u' V^3\ge0$.
\end{pf}

For the remainder of the proof of Theorem~\ref{thm:3}, let $e$ denote
the three-vector of
ones, so (\ref{eq:52}) is equivalently expressed as $\theta=-e$, and
we can represent $X$ in (\ref{11}) as
%
\begin{equation}
\label{eq:56}
X(t)=X(0) + B(t) - e t,\qquad    t \ge0,
\end{equation}
where $B$ is a driftless Brownian motion with nonsingular covariance matrix
and $B(0) = 0$. Also, we choose a starting state $x=X(0) = Z(0)$
that satisfies
%
\begin{equation}
\label{eq:57}
Z_{1}(0) \ge0,\qquad    Z_{2}(0) = 0\quad   \mbox{and}\quad  Z_{3}(0) > 0.
\end{equation}
In this section, because the initial state is fixed, we write
$\mathbb{E}(\cdot)$ rather than $\mathbb{E}_{x}(\cdot)$ to signify
the expectation
operator associated with the probability measure $\mathbb{P}_{x}$ (see
Appendix \hyperref[sec:srbm]{A}). Also, when we speak of stopping times and
martingales, the relevant filtration is the one specified in
Appendix~\hyperref[sec:srbm]{A}.

Let $u> 0$ be chosen to satisfy
${u}'V\ge0$, as in Lemma \ref{lem:51}, and further normalized so that
${u}'e=1$. It is immediate from the definition of $V$ that
${u}'V={u}'R-{e}'$, and thus one has the following:
%
\begin{equation}
\label{eq:58}
{u}'R\ge{e}' .
\end{equation}
Now define $\xi(t)={u}'Z(t)$, $t \ge$ 0. From (\ref{11}), (\ref
{eq:56}) and
(\ref{eq:58}), one has
\begin{eqnarray}
\label{eq:59}
\xi(t)-\xi(0) &=& {u}'B(t)-{u}'e t+{u}'RY(t)\nonumber\\[-8pt]\\[-8pt]
&\ge&
{u}'B(t)-t+{e}'Y(t)\qquad  \mbox{for } t \ge0.\nonumber
\end{eqnarray}
Next, let
\begin{eqnarray*}
\tau_{1}&=& \inf\{t > 0\dvtx  Z_{3}(t) = 0\},\qquad   \tau_{2} = \inf\{t
> \tau_{1}\dvtx  Z_{1}(t) = 0\},\\
 \tau_{3} &=& \inf\{t > \tau
_{2}\dvtx  Z_{2}(t) = 0\}
\end{eqnarray*}
and so forth. (These stopping times are analogous to the points in
time at which the piecewise linear fluid path in Figure \ref{fig2}
changes direction.) The~crucial observation is the following:
$Z_{3}(\cdot) > 0 $ over the interval $[0,\tau_{1})$,
$Z_{1}(\cdot) > 0 $ over $[\tau_{1},\tau_{2})$, \mbox{$Z_{2}(\cdot
) > 0$} over $[\tau_{2},\tau_{3})$ and so forth. Thus
$Y_{3}(\cdot)$ does not increase over $[0,\tau_{1})$, \mbox{$Y_{1}(\cdot
)$} does not increase over $[\tau_{1},\tau_{2})$, \mbox{$Y_{2}(\cdot)$}
does not increase over $[\tau_{2},\tau_{3})$ and so forth.

From
(\ref{11}) and (\ref{eq:56}), we then have the following relationships:
\begin{eqnarray}
\label{eq:510}
Z_{2}(t)&=&B_{2}(t)- t+Y_{2}(t)\nonumber \\[-8pt]\\[-8pt]
&& {} + R_{21} Y_{1}(t) ,\qquad   0\le t \le
\tau_{1},\nonumber\\
\label{eq:511}
Z_{3}(t) &= & [B_{3}(t) - B_{3}(\tau_{1})] - (t - \tau_{1})
+ [Y_{3}(t)- Y_{3}(\tau_{1})] \nonumber\\[-8pt]\\[-8pt]
 && {} + R_{32}[Y_{2}(t)- Y_{2}(\tau
_{1})] ,\qquad   \tau_{1} \le t \le\tau_{2}, \nonumber\\
\label{eq:512}
Z_{1}(t) & = & [B_{1}(t) - B_{1}(\tau_{2})] - (t - \tau_{2})
+ [Y_{1}(t) - Y_{1}(\tau_{1})] \nonumber\\[-8pt]\\[-8pt]
 && {} + R_{13}[Y_{3}(t)- Y_{3}(\tau
_{2})], \qquad  \tau_{2} \le t \le\tau_{3}. \nonumber
\end{eqnarray}
There exist analogous representations for
$Z_2$ over the time interval $[\tau_3, \tau_4]$, for $Z_3$ over
$[\tau_4, \tau_5]$, for $Z_1$ over $[\tau_5, \tau_6]$, and so on.
Now (\ref{eq:510}) gives
%
\begin{equation}
\label{eq:513}
Y_{2}(t) = t+Z_{2}(t)- B_{2}(t) - R_{21}Y_{1}(t)\qquad   \mbox{for } 0
\le
t \le\tau_{1}.
\end{equation}
Because $Y_{3} \equiv0$ on $[0,\tau_{1})$, we can substitute (\ref{eq:513})
into (\ref{eq:59}) to obtain the following:
\begin{eqnarray}
\label{eq:514}
\hspace*{15pt} \xi(t)-\xi(0)&\ge&{u}'B(t)-t+Y_{1}(t)\nonumber\\[-8pt]\\[-8pt]
\hspace*{15pt} &&{} +
[t+Z_{2}(t)- B_{2}(t)- R_{21}Y_{1}(t)] \qquad  \mbox{for 0 $\le t
\le\tau_{1}$}.\nonumber
\end{eqnarray}
From the definition of $V$ and (\ref{eq:54}), we have $1 - R_{21}=
b_1 >0$, so (\ref{eq:514}) can be rewritten
%
\begin{equation}
\label{eq:515}
\xi(t)-\xi(0)\ge M(t)+A(t)\qquad   \mbox{for } 0 \le t \le\tau_{1},
\end{equation}
where
%
\begin{eqnarray}
\label{eq:516}
 M(t)&=&{u}'B(t)- B_2 (t)\qquad   \mbox{for } 0 \le t \le\tau_{1},
\\
 A(t)&=&Z_2 (t)+ b_1 Y_1 (t)\qquad   \mbox{for } 0 \le t \le\tau
_{1}. \label{eq:517}
\end{eqnarray}
Defining $\tau= \lim\tau_{n}$, we now extend the definition
(\ref{eq:516}) to
all $t\in[0,\tau)$ as follows:
\begin{eqnarray}
\label{eq:518}
M(t)& =&M(\tau_{1})+{u}'[B(t)-B( \tau_1
)]\nonumber\\[-8pt]\\[-8pt]
&&{}- {[B}_3 (t)- B_3 {(\tau
}_1 )]\qquad  \mbox{for } \tau_{1} \le t \le\tau_{2},\nonumber
\label{eq:519} \\
 M(t)&=&M(\tau_{2})+{u}'[B(t)-B( \tau_2
)]\nonumber\\[-8pt]\\[-8pt]
&&{}- {[B}_1 (t)- B_1 {(\tau
}_2 )] \qquad  \mbox{for } \tau_{2} \le t \le\tau_{3}\nonumber
\end{eqnarray}
and so forth. Finally, on $\{\tau< \infty\}$, we set $M(t) =
M(\tau)$ for all $t \ge\tau$. Then $M = \{M(t),  t \ge0\}$ is a
continuous martingale whose quadratic variation $\langle
M, M\rangle(\cdot)$ satisfies
%
\begin{equation}
\label{eq:520}
\langle M,M\rangle(t)-\langle M,M\rangle(s) \le\gamma(t -
s)\qquad   \mbox{for }0 < s < t < \infty,
\end{equation}
where $0 < \gamma< \infty$.
Also, we extend (\ref{eq:517}) to all $t\in[0,\tau)$ via
%
\begin{eqnarray}
\label{eq:521}
 A(t)&=&A(\tau_{1})+Z_{3}(t)+b_{2}[Y_{2}(t) -Y_{2}(\tau
_{1})]\qquad   \mbox{for } \tau_{1} \le t \le\tau_{2}, \\
 A(t)&=&A(\tau_{2})+Z_{1}(t)+b_{3}[Y_{3}(t) -Y_{3}(\tau
_{2})]\qquad   \mbox{for } \tau_{2} \le t \le\tau_{3} \label{eq:522}
\end{eqnarray}
and so forth. Thus the process $A = \{A(t), 0 \le t < \tau\}$ is
nonnegative and continuous.

\begin{lemma}\label{lem:52}
$\xi(t)-\xi(0) \ge M(t)+A(t)$ for all $t\in
[0,\tau)$.
\end{lemma}

\begin{pf}
It has already been shown in (\ref{eq:515}) that this inequality is
valid for $0 \le
t \le\tau_{1}$. In exactly the same way, but using (\ref{eq:511}) instead
of (\ref{eq:510}), one obtains
%
\begin{equation}
\label{eq:523}
\qquad \quad \xi(t)-\xi(\tau_{1}) = [M(t)-M(\tau_{1})] + [A(t)-
A(\tau_{1})] \qquad  \mbox{for } \tau_{1} \le t \le\tau_{2},
\end{equation}
so the desired inequality holds for $0 \le t \le\tau_{2}$.
Continuing in this way, the desired inequality is established for $0\le
t < \tau$.
\end{pf}

To complete the proof of Theorem \ref{thm:3}, let $T = \inf\{t > 0\dvtx
\xi(t) = \epsilon\}$ and let $\sigma= \inf\{t > 0\dvtx  \xi(0) +
M(t) = \epsilon\}$, where $0< \epsilon< \xi(0)$. From Lemma
\ref{lem:52}, the nonnegativity of $A(\cdot)$, and the fact that
$\xi(\tau) = 0$ on $\{\tau< \infty\}$, we have the following
inequalities: $0< \sigma\le T \le\tau$.
Thus it suffices to prove that $\mathbb{E}(\sigma) =\infty$, which
can be
shown by
essentially the same argument that applies when $M$ is an ordinary
(driftless) Brownian motion. That is, we first let $\sigma(b) = \inf\{
t > 0\dvtx
\xi(0)+ M(t) = \epsilon\mbox{ or } \xi(0) + M(t) = b\}$, where $b
>\xi(0)$. Because both $M$ and $M^2-\langle M,M\rangle$
are martingales and (\ref{eq:520}) holds, one has $0
< \mathbb{E}[\sigma(b)] <\infty$,
$\mathbb{E}[M(\sigma(b))] =0$ and
$\mathbb{E}[M^2(\sigma(b))]=\mathbb{E}[\langle M,M\rangle(\sigma
(b))]$. It follows by
the optional sampling theorem that
%
\begin{eqnarray}\label{eq:524}
\mathbb{E}[\langle M,M\rangle(\sigma(b))] &=&
\mathbb{E}[M^2(\sigma(b))]\nonumber\\
&=&
\bigl(b-\xi(0)\bigr)^2\frac{\xi(0)-\epsilon}{b-\epsilon}+\bigl(\xi(0)-\epsilon
\bigr)^2\frac
{b-\xi(0)}{b-\epsilon}\\
&= & \bigl(b-\xi(0)\bigr)\bigl(\xi(0)-\epsilon\bigr).\nonumber
\end{eqnarray}
The~left-hand side of (\ref{eq:524}) is $\le\gamma E[\sigma(b)]$ by
(\ref{eq:520}), the right-hand side $\uparrow\infty$ as $b\uparrow
\infty$, and obviously $\sigma\ge\sigma(b)$ for all $b >
\xi(0)$. Thus $\mathbb{E}(\sigma) =\infty$, and the proof
of Theorem \ref{thm:3} is complete.


\section{Categories of divergent LCP solutions}
\label{sec:divergent}
Our goal in the remainder of the paper is to prove Theorem~\ref{thm:4}.
We continue to assume the canonical problem format in which $R$
satisfies (\ref{eq:new1.5}) and $\theta$ satisfies (\ref{eq:51}). In
the following lemma and later, the term ``LCP solution'' is used to
mean a solution $(u,v)$ of the linear complementarity problem
(\ref{eq:LCP1})--(\ref{eq:LCP3}).

\begin{lemma} \label{lem:3}
If (\ref{15}) holds, then \textup{(a)} $\theta\ge0$ is not possible, and \textup{(b)} there
exists no LCP solution
$(u,v)$ with $v>0$.
\end{lemma}

\begin{pf}
Because $R$ is completely $\mathcaligr{S}$ by assumption, its transpose is
also completely $\mathcaligr{S}$; cf. Proposition 1.1 of Dai and Williams
(\citeyear{daiwil95}). Thus there exists a vector $a > 0$ such that $a'R > 0$. Now
(\ref{15}) says that $\theta+ Ry = 0$ for $y > 0$. Multiplying
both sides of the equation by $a'$ and rearranging terms, one has
$a'\theta=- a'R y< 0$, which implies conclusion (a). Also, if
$(u,v)$ is a LCP solution with
$v>0$, one has from
(\ref{eq:LCP3}) and (\ref{eq:LCP2}) that $u=0$ and $v
=\theta$, which contradicts conclusion (a). This implies conclusion (b).
\end{pf}

We now define five nonoverlapping categories of divergent LCP solutions.
Immediately
after each category is defined, we shall exhibit a pair $(R, \theta)$
which admits a LCP solution $(u,v)$ in that category, or else direct
the reader to a proposition that shows the category to be empty. Readers
may verify that the reflection matrix $R$ appearing in each of our
examples is completely $\mathcaligr{S}$.
Also, defining $u^*=-R^{-1}\theta$ as in (\ref{eq:ustar}), we shall
display the vector $u^*$ for each example, showing that $u^*>0$ and
hence (\ref{15}) is satisfied.

\renewcommand{\thecategory}{\Roman{category}}
\setcounter{category}{0}
\begin{category}\label{cat1}
Exactly two components of $v$ are positive, and
the complementary component of $u$ is positive.
The~following is such an example:
\begin{eqnarray*}
R&=&
\pmatrix{
1 & 1/3 & 1/3 \cr
2 & 1 & -1/2 \cr
2 & -1/2 & 1},\qquad
 \theta=
\pmatrix{
-1 \cr -1 \cr -1},\\
u^* &=&
\pmatrix{
1/5 \cr 6/5 \cr 6/5},\qquad
  u=
\pmatrix{
1 \cr 0 \cr 0},\qquad
  v=
\pmatrix{
0 \cr 1 \cr 1}.
\end{eqnarray*}
\end{category}

\begin{category}\label{cat2}
 Exactly one component of $v$ is
positive, $\det(\hat R)>0$,
and the two complementary components of $u$ are not both zero, where
$\hat R$ is the $2\times2$ principal submatrix of $R$
corresponding to the two zero components of $v$. Such an example is
given by
\begin{eqnarray*}
R&=&
\pmatrix{
1 & 1 & 1/2 \cr
-2 & 1 & 0 \cr
3 & 0 & 1},\qquad
 \theta=
\pmatrix{
-1 \cr 1 \cr -1},\\
u^*&=&
\pmatrix{
1 \cr 1 \cr2},\qquad
  u=
\pmatrix{
2/3 \cr 1/3 \cr 0},\qquad
  v=
\pmatrix{
0 \cr 0 \cr 1}.
\end{eqnarray*}
\end{category}

Here, the two complementary components of $u$ are both
positive. In the following example, which also falls in Category \ref{cat2},
just one of them is positive:
\begin{eqnarray*}
R&=&
\pmatrix{
1 & 1/2 & 3 \cr
1 & 1 & 2 \cr
2 & 1 & 1},\qquad
 \theta=
\pmatrix{
-1 \cr -1 \cr -1},\\
u^* &=&
\pmatrix{
1/5 \cr 2/5 \cr 1/5},\qquad
  u=
\pmatrix{
1 \cr 0 \cr 0},\qquad
  v=
\pmatrix{
0 \cr 0 \cr 1}.
\end{eqnarray*}

\begin{category}\label{cat3}
 Exactly one component of $v$ is positive,
$\det(\hat R)=0$, and the two complementary components of $u$ are
not both zero. In Lemma~\ref{lem:new6}, it will be shown that no
such LCP solutions exist if (\ref{15}) holds.
\end{category}

\begin{category}\label{cat4}
 Exactly one component of $v$ is positive,
$\det(\hat R)<0$, and the two complementary components of $u$ are
both positive. Such an example is given by
\begin{eqnarray*}
R&=&
\pmatrix{
1 & 11/10 & 2 \cr
2 & 1 & 0 \cr
0 & 2 & 1},\qquad
  \theta=
\pmatrix{
-1 \cr -1 \cr -1},\\
u^*&=&
\pmatrix{
19/68 \cr
15/34 \cr
2/17},\qquad
  u=
\pmatrix{
1/12 \cr 5/6 \cr 0},\qquad
  v=
\pmatrix{
0 \cr 0 \cr 2/3}.
\end{eqnarray*}
\end{category}

It will be shown in Lemma \ref{lem:caseIIIB} that if there
exists a LCP solution in Category~\ref{cat4}, under our restrictions on
$R$ and $\theta$, there also exists a solution in Category~\ref{cat1} or
Category \ref{cat2} (or both). For the example above, a second LCP solution is
$(\hat u, \hat v)$, where $\hat u = (0,1,0)'$ and $\hat v =
(1/10,0,1)'$; this second solution lies in Category \ref{cat1}.

\begin{category}\label{cat5}
Exactly one component of $v$ is
positive, $\det(\hat R)<0$, and exactly one of the two complementary
components of $u$ is positive. Such an example is given by
\begin{eqnarray*}
R&=&
\pmatrix{
1 & 1 & -2/5 \cr
2 & 1 & -6/5 \cr
-2 & -1/10 & 1},\qquad
  \theta=
\pmatrix{
-1 \cr -1 \cr 1},\\
u^*&=&
\pmatrix{
9/8\cr 5/14 \cr 45/28},\qquad
  u=
\pmatrix{
0 \cr 1 \cr 0},\qquad
  v=
\pmatrix{
0 \cr 0 \cr 0.9}.
\end{eqnarray*}
\end{category}

\begin{lemma}\label{lem:classificaitionOverview}
Suppose that (\ref{15}) holds and that $(u,v)$ is a divergent LCP
solution. Then $(u,v)$ belongs to one of the five categories defined
immediately above.
\end{lemma}

\begin{pf}
Let $m$ and $n$ denote the number of positive components in $u$ and
$v$, respectively; the complementarity condition (\ref{eq:LCP3})
implies that $m + n \le3$. Lemma \ref{lem:3} shows that $n < 3$
(i.e., $v>0$ cannot hold); also $n > 0$, because $(u,v)$ is a
divergent LCP solution by assumption. Thus, either $n = 1$ or $n =
2$. Moreover, it is not possible that $m=0$, or equivalently $u =
0$, because then (\ref{eq:LCP1}) and (\ref{eq:LCP2}) would imply
$\theta= v \ge0$, which contradicts Lemma \ref{lem:3}. So the
only remaining possibilities are $(m,n) = (1,2)$, $(m,n) = (2,1)$
and $(m,n) = (1,1)$. Category \ref{cat1} is precisely the case where $(m,n) =
(1,2)$, and Categories \ref{cat2} through \ref{cat5} together cover the cases where
$(m,n) = (2,1)$ and $(m,n) = (1,1)$.
\end{pf}

It will be shown in Section \ref{sec:lemmas} that $Z$ cannot be
positive recurrent if there exists a LCP solution in Category \ref{cat1}, Category \ref{cat2} or Category \ref{cat5}. Lemma~\ref{lem:caseIIIB} in this
section will show
that the existence of a LCP solution in Category \ref{cat4} implies the
existence of a LCP solution in either Category \ref{cat1} or Category \ref{cat2}.
Lemma \ref{lem:new6} at the end of this section
will show that LCP solutions in Category \ref{cat3} cannot occur when
\eqref{15} holds. In combination with Lemma
\ref{lem:classificaitionOverview} above, these results obviously imply
Theorem \ref{thm:4}.

We now state and prove Lemma~\ref{lem:theta}, which we need in order
to prove Lemma~\ref{lem:caseIIIB}. Our scaling convention
(\ref{eq:new1.5}) specifies that $R$ has ones on the diagonal, so we
can write
%
\begin{equation}
\label{eq:Rspecial}
R=
\pmatrix{
1 & a' & c \cr
a & 1 & c' \cr
b & b' & 1}
\end{equation}
for some constants $a, a'$, $b, b'$, $c$ and $c'$.


\begin{lemma}\label{lem:theta}
Assume that there does not exist a LCP solution in
Category \textup{\ref{cat1}}, and that there is a divergent LCP solution $(u, v)$ with
$u_1>0$, $u_2> 0$, $u_3=0$, $v_1=v_2=0$ and $v_3>0$.
Let $R$ be as in (\ref{eq:Rspecial}) and assume that the principal
submatrix $\hat R$ corresponding to the zero components of $v$
satisfies $\det(\hat R)<0$. Then $\theta=(-1, -1, 1)'$ and $a,
a'>1$.
\end{lemma}

\begin{pf}
Because $(u,v)$ is a solution of the LCP
(\ref{eq:LCP1})--(\ref{eq:LCP3}), one has
%
\begin{equation}
\label{eq:specialLCP}
\pmatrix{
1 & a' & c \cr
a & 1 & c' \cr
b & b' & 1}
\pmatrix{
u_1 \cr u_2 \cr 0}
=
\pmatrix{
-\theta_1 \cr -\theta_2 \cr -\theta_3+v_3}.
\end{equation}
%
Because $v_1=v_2=0$ and $v_3>0$,
%
\begin{equation}
\label{eq:RhatDef}
\hat R =
\pmatrix{
1 & a' \cr a & 1}.
\end{equation}
Setting $\hat{u}=( u _1 , u _2 )'$ and
$\hat{\theta}=( \theta_1, \theta_2 )'$,
we have from (\ref{eq:specialLCP}) that
%
\begin{equation}
\label{eq:l3}
\hat{R} \hat{u}=-\hat{\theta}.
\end{equation}
Because $\hat{R}$ is an
$\mathcaligr{S}$-matrix with negative determinant,
%
\begin{equation}
\label{eq:aPBigerOne}
a, a'>0\quad   \mbox{and}\quad    aa'>1.
\end{equation}
Because $u_1>0$ and $u_2> 0$ by hypothesis, it is immediate from
(\ref{eq:l3}) and (\ref{eq:aPBigerOne}) that both components of $\hat
\theta$ are negative, so our
canonical rescaling gives $\hat{\theta}=(-1,-1)'$.
Thus, either $\theta=(-1, -1, -1)'$, $(-1, -1, 0)'$ or $(-1, -1, 1)'$
must hold.
From
(\ref{eq:aPBigerOne}), (\ref{eq:l3}) and $\hat{\theta}=(-1,-1)'$,
it follows that
%
\begin{equation}
\label{eq:aBiggerOne}
a,{a}'>1.
\end{equation}

We will show that $\theta=(-1, -1, 1)'$ by excluding the other two
cases. Suppose first that $\theta=(-1,-1,-1)'$. Then
(\ref{eq:specialLCP}) becomes
%
\begin{equation}
\label{eq:specialLCP1}
\pmatrix{
1 & a' & c \cr
a & 1 & c' \cr
b & b' & 1}
\pmatrix{
u_1 \cr u_2 \cr 0}
=
\pmatrix{
1\cr 1 \cr 1+v_3}.
\end{equation}
It must be true that $b,{b}'\le1$; otherwise there would be a
solution of the LCP that falls into Category \ref{cat1}. For example, if $b>1$
one then has a divergent LCP solution $(\bar u, \bar v)$ with $\bar
u=(1, 0, 0)'$, $\bar v_1=0$, $\bar v_2=a-1>0$ and $\bar{v}_3=b-1>0$.
However, one cannot have $a, a'\ge1$, $b,{b}'\le1$ and
(\ref{eq:specialLCP1}) holding simultaneously, which gives a
contradiction.

Next suppose that $\theta=(-1,-1,0)'$. Here, (\ref{eq:specialLCP1})
holds with $v_{3 }$ in place of $1+v_{3 }$ on the right-hand side.
We must have $b,{b}'\le0$, for the same reason as before.
This results in a contradiction, and the only remaining possibility
under our canonical
rescaling is $\theta=(-1,-1,1)'$.
\end{pf}

\begin{lemma} \label{lem:caseIIIB}
If there exists a LCP solution in
Category \textup{\ref{cat4}}, then there also exists a solution in Category \textup{\ref{cat1}} or
Category \textup{\ref{cat2}} (or both).
\end{lemma}

\begin{pf}
Denoting by $(u, v)$ a solution in Category \ref{cat4}, we assume that no
solution in Category~\ref{cat1} exists. It will then suffice to prove that a
solution in Category~\ref{cat2} exists. By permuting the indices, we can
assume that $u_1>0$,
$u_2>0$, $u_3=0$, $v_1=v_2=0$ and $v_3>0$.

We use the
notation in (\ref{eq:Rspecial}).
By Lemma~\ref{lem:theta},
one has $a, a^\prime> 1$ and $\theta= (-1, -1,
1)^\prime$. We shall assume that $c^\prime\ge
c$, then construct a LCP solution $(\tilde{u},\tilde{v})$ that
falls into Category \ref{cat2}, with $\tilde{u}_1>0$ ,
$\tilde{u}_3 > 0$, $\tilde{v}_2>0$; in exactly the same way, if $c
\ge c^\prime$,
one can construct a LCP solution $(\tilde{u},\tilde{v})$
that falls into Category \ref{cat2}, with ${\tilde{v}}_1>0$ , ${\tilde
{u}}_2> 0$ , ${\tilde{u}}_3 >0$.

We first observe that $b, b^\prime\le-1$. Otherwise,
contrary to the assumption imposed in the first paragraph of the
proof, there would exist a LCP solution in Category~\ref{cat1}.
For example, if $b > -1$, then there is a divergent LCP solution
$(\bar{u},\bar{v})$ with $\bar u = (1,0,0)',$ $\bar{v}_1=0$, $\bar
v_2=a-1>0$, and $\bar v_3=b+1>0$.


The~$2\times2$
submatrix of $R$ that is relevant to our construction is
\[
\tilde R=
\pmatrix{
1&c \cr b &1}.
\]
Because $\tilde{R}$ is an $\mathcaligr{S}$-matrix and $b< 0$, we
know that $bc < 1$, hence $\det(\tilde R)>0$, and because $b\le-1$
we also know that
$c > -1$. Letting $\gamma=(\gamma_1
, \gamma_2 )'$ be the two-vector satisfying $\tilde
{R}\gamma=(1,-1)'$, one has
\[
\gamma=\frac{1}{1-bc}
\pmatrix{
1 & -c \cr -b & 1
}
\pmatrix{
1 \cr -1
}
=\frac{1}{1-bc}
\pmatrix{
1+c\cr -1 -b
}
.
\]
Defining $\tilde{u}=( \gamma_1 ,0, \gamma
_2 )'$ and $\tilde{v}=\theta+R\tilde{u}$,
it follows that
$ {\tilde{v}}_{1} = {\tilde
{v}}_{ 3} =0$. Comparing the
first and second rows of $R$ term by term, and noting that the first two
components of $\theta$ are identical, one sees that
%
\begin{equation}
\label{eq2}
{\tilde{v}}_2 - {\tilde{v}}_1 =\frac{1}{1-bc}[
(a-1)(1+c)-(c'-c)(1+b) ].
\end{equation}
Because of the inequalities $a>1$, $c > -1$, $c' \ge{c}$ and $b\le
-1$, the quantity inside the square brackets in (\ref{eq2}) is
positive. Thus $ {\tilde{v}}_2 >0$, and hence $(\tilde
{u},\tilde{v})$ is a LCP solution in Category \ref{cat2}.
\end{pf}

\begin{lemma} \label{lem:new6} If (\ref{15}) holds, then there cannot
exist a LCP solution in Category~\textup{\ref{cat3}}.
\end{lemma}

\begin{pf}
Arguing as in the proof of Lemma \ref{lem:caseIIIB}, we assume the
existence of a LCP solution $(u,v)$ in Category \ref{cat3}. By permuting
the indices, we can assume that $u_1 > 0$, $u_2\ge0$, $u_3 = 0$,
$v_1 = v_2 = 0$ and $v_3 > 0$. We use the notation (\ref{eq:Rspecial})
and define
$\hat R$ by (\ref{eq:RhatDef}). A minor variation of the first
paragraph in the proof of Lemma~\ref{lem:theta} shows for the
current case that both $\theta_1$ and $\theta_2$ are negative, and
so $\theta_1=\theta_2=-1$ with our scaling convention. One then has
$a = a' = 1$ in (\ref{eq:Rspecial}), because $\det(\hat R) = 0$ by
assumption. By (\ref{15}), $\theta+ Ru^* = 0$ for some $u^* > 0$,
from which it follows that $c = c'$ in (\ref{eq:Rspecial}). That is,
the first two rows of $R$ are identical, whereas (\ref{15}) includes
the requirement that $R$ be nonsingular.
\end{pf}

\section{Proof of Theorem \protect\ref{thm:4}}
\label{sec:lemmas}

As we explained immediately after the proof of Lemma
\ref{lem:classificaitionOverview} in Section \ref{sec:divergent}, the
proofs of Lemmas~\ref{lem:caseI}, \ref{lem:caseII} and
\ref{lem:caseV} in this section will complete the proof of
Theorem~\ref{thm:4}. In Lemmas~\ref{lem:caseI} and \ref{lem:caseII},
we actually prove that $Z$ is transient, which is stronger than
we require for Theorem~\ref{thm:4}. The~SRBM $Z$
is said to be \textit{transient} if there exists an open ball $C$
centered at the origin such that $\mathbb{P}\{\tau_C=\infty\}>0$ for some
initial state $Z(0)=x\in\mathbb{R}^3_+$ that is outside of the ball, where
$\tau_C=\inf\{t\ge0\dvtx  Z(t)\in C\}$. Clearly, when $Z$ is transient, it
is not positive recurrent. In this section, we continue to assume the
canonical problem format in which $R$ satisfies (\ref{eq:new1.5}) and
$\theta$ satisfies (\ref{eq:51}).

\begin{lemma} \label{lem:caseI}
If there
is a LCP solution $(u, v)$ in Category \textup{\ref{cat1}}, then $Z$ is transient.
\end{lemma}

\begin{pf}
Without loss of generality, we assume that
%
\begin{equation}
\qquad u_1>0,\qquad    u_2=0,\qquad    u_3=0, \qquad   v_1=0,\qquad    v_2>0, \qquad  v_3>0.
\end{equation}
Because $v=\theta+Ru$, with $R$ as in (\ref{eq:Rspecial}),
$\theta_1<0$, and so,
by our scaling convention, $\theta_1=-1$.
It follows from this that
%
\begin{equation}
\label{eq:CategoryIv3}
 u_1=1,\qquad
v_2=\theta_2+ a> 0\quad   \mbox{and} \quad   v_3=\theta_3+b>0.
\end{equation}

One can write (\ref{11}) as
%
\begin{eqnarray}
 Z_1(t) &=& Z_1(0)+ \theta_1 t + B_1(t) + Y_1(t) + a' Y_2(t)+c
Y_3(t), \label{eq:Z1}\\
 Z_2(t)& =& Z_2(0) + \theta_2 t + B_2(t) + a Y_1(t)+ Y_2(t) +c'
Y_3(t), \label{eq:Z2} \\
 Z_3(t)& =& Z_3(0) + \theta_3 t + B_3(t) + b
Y_1(t)+b'Y_2(t) + Y_3(t) \label{eq:Z3}
\end{eqnarray}
for $t\ge0$,
where $B=\{B(t), t\ge0\}$ is the three-dimensional driftless Brownian
motion with covariance matrix $\Gamma$. Assume $Z(0)=(0, N,
N)'$ for some
constant $N>1$ and set $\tau= \inf\{t\ge0\dvtx
Z_2(t)=1 \mbox{ or } Z_3(t)=1\}$. We will show that
$\mathbb{P}\{\tau=\infty\}>0$ for sufficiently large $N$, and thus
$Z$ is
transient.
Because $\theta_1=-1$ and $Y_2(t)=Y_3(t)=0$ for $t\in[0,
\tau)$, one has, for $t< \tau$,
\begin{eqnarray*}
 Z_1(t) &=& -t + B_1(t) + Y_1(t), \\
 Z_2(t)& =& N + \theta_2 t + B_2(t) + a Y_1(t), \\
 Z_3(t)& =& N+ \theta_3 t + B_3(t) + b Y_1(t).
\end{eqnarray*}
For $t\ge0$, let $\hat Y_1(t) = \sup_{0\le s\le t}(-s+B_1(s))^-$, and
set
%
\begin{eqnarray}
 \hat Z_1(t)& =& -t + B_1(t) + \hat Y_1(t), \label{eq:hZ1}\\
 \hat Z_2(t) &=& N + \theta_2 t + B_2(t) + a \hat Y_1(t)
, \label{eq:hZ2} \\
 \hat Z_3(t) &=& N + \theta_3 t + B_3(t) + b \hat Y_1(t)
\label{eq:hZ3}
\end{eqnarray}
for $t\ge0$.
Clearly, $Z(t)=\hat Z(t)$ for $t\in[0, \tau]$. In particular, $\tau
=\hat\tau$,
where $\hat\tau= \inf\{t\ge0\dvtx  \hat Z_2(t)=1 \mbox{ or } \hat
Z_3(t)=1\}$. To show $Z$ is transient, it suffices to prove that, for
sufficiently large $N$,
$\mathbb{P}\{\hat\tau=\infty\}>0$.

By the functional
strong law of large numbers (FSLLN) for
a driftless Brownian motion, one has
\[
\lim_{t\to\infty} \frac{1}{t}\sup_{0\le s \le1}|B(ts)|=0\qquad
\mbox{almost surely}.
\]
This implies that
\begin{eqnarray*}
\lim_{t\to\infty} t^{-1}\hat Y_1(t)
&= & \lim_{t\to\infty} t^{-1} \sup_{0\le s\le t } \bigl(-s+B_1(s)\bigr)^{-} \\
&= & \lim_{t\to\infty} \sup_{0\le s\le1 } \bigl(-s+t^{-1}B_1(ts)\bigr)^- \\
&= & \sup_{0\le s \le1} (-s)^{-} = 1\qquad
\mbox{almost surely.}
\end{eqnarray*}
%
Therefore, by (\ref{eq:CategoryIv3}), (\ref{eq:hZ2}) and
(\ref{eq:hZ3}), one has
$\lim_{t\to\infty} t^{-1} \hat Z_2(t)=v_2>0$ and\vspace*{1pt}
$\lim_{t\to\infty} t^{-1} \hat Z_3(t)=v_3>0$ almost surely.
Consequently, there exists a constant $T>0$ such that $\mathbb
{P}(A)\ge
3/4$
for all $N\ge1$, where
%
\begin{equation}
\label{eq:78a}
A =\{ \hat Z_2(t)>1 \mbox{ and } \hat Z_3(t)>1 \mbox{ for all $t\ge
T$}\}.
\end{equation}
One can choose $N$ large enough so that $\mathbb{P}(B)\ge3/4$, where
%
\begin{equation}
\label{eq:78b}
B=\{\hat Z_2(t)>1 \mbox{ and } \hat Z_3(t)>1 \mbox{ for all
$t\in[0, T]$}\}.
\end{equation}
Because $A\cap B \subset\{ \hat\tau=\infty\}$, $\mathbb{P}\{\hat
\tau=\infty\}\ge1/2>0$, as desired.
\end{pf}

\begin{lemma} \label{lem:caseII} If there
is an LCP solution $(u, v)$ in Category \textup{\ref{cat2}}, then $Z$ is
transient.
\end{lemma}

\begin{pf}
Without loss of generality, we assume that
%
\begin{equation}
\label{eq:uvIIIA}
\qquad u_1>0,\qquad    u_2\ge0,\qquad    u_3=0,\qquad    v_1=0, \qquad   v_2=0, \qquad  v_3>0.
\end{equation}
Assume $R$ is given by (\ref{eq:Rspecial}), and let $\hat R$
be the $2\times2$ principal submatrix of $R$ given by
(\ref{eq:RhatDef}). By assumption, $\det(\hat R)>0$. One can check
that conditions
(\ref{eq:LCP1})--(\ref{eq:LCP3}) and (\ref{eq:uvIIIA}) imply that
%
\begin{eqnarray}
\hat R^{-1}
\pmatrix{
\theta_1 \cr
\theta_2
}
& \le&0, \label{eq:LCPv2}\\
 v_3 &=& \theta_3 - (b, b')\hat R^{-1}
\pmatrix{
\theta_1 \cr \theta_2
}
>0. \label{eq:LCPv3}
\end{eqnarray}

Let $Z(0)=(0, 0, N)'$ for some
constant $N>1$ and set $\tau=
\inf\{t\ge0\dvtx
Z_3(t)=1\}$.
We will show that for sufficiently large $N$,
$\mathbb{P}\{\tau=\infty\}>0$, and thus $Z$ is transient.

On $\{t<\tau\}$, one has $Z_3(t)>0$ and thus $Y_3(t)=0$.
Because the SRBM $Z$ satisfies Equation (\ref{11}),
on $\{t<\tau\}$,
\begin{eqnarray}
\label{eq:RtildeSRBM}
Z(t) &=& Z(0) + \theta t + B(t) + R
\pmatrix{
Y_1(t) \cr Y_2(t) \cr 0
}\nonumber\\[-8pt]\\[-8pt]
&= &Z(0) + \theta t + B(t) +\tilde R
\pmatrix{
Y_1(t) \cr Y_2(t) \cr Y_3(t)
}
,\nonumber
\end{eqnarray}
where
\[
\tilde R =
\pmatrix{
1 & a' & 0 \cr
a & 1 & 0 \cr
b & b' & 1
}
.
\]
One can check that because $R$ is completely $\mathcaligr{S}$, so is $\tilde
R$. It therefore follows from \citet{taywil93} that there exists a
SRBM $\tilde Z$ associated with the data $(\mathbb{R}^3_+,\theta,
\Gamma,
\tilde R)$ that starts from $Z(0)$. Following Definition
\ref{def:srbm} in Appendix~\hyperref[sec:srbm]{A}, the three-dimensional
process $\tilde Z$, together
with the corresponding processes $\tilde B$ and $\tilde Y$,\vspace*{1pt} is defined on
some filtered probability space $(\tilde\Omega, \{ \tilde\mathcaligr{F}_t\},\tilde\mathbb{P})$; $\tilde B$, $\tilde Y$ and $\tilde Z$ are
adapted to $\{\tilde\mathcaligr{F}_t\}$, and $\tilde\mathbb{P}$-almost surely
satisfy (\ref{11})--(\ref{14}); $\tilde B$ is a driftless Brownian
motion with covariance matrix $\Gamma$, and $\tilde B$ is an $\{\tilde
\mathcaligr{F}_t\}$-martingale. Furthermore, from \citet{taywil93}, the
distribution of $\tilde Z$ is unique. Because of
(\ref{eq:RtildeSRBM}), $B$, $Y$ and $Z$ also satisfy
(\ref{11})--(\ref{14}) on $\{t< \tau\}$, with the same data
$(\mathbb{R}^3_+,
\theta, \Gamma, \tilde R)$, and so $\tau=\tilde\tau$ in
distribution, where
\[
\tilde\tau= \inf\{t\ge0\dvtx  \tilde Z_3(t)=1\}.
\]

We now show that for sufficiently large $N$,
%
\begin{equation}
\label{eq:tauhat}
\tilde\mathbb{P}\{ \tilde\tau=\infty\}>0,
\end{equation}
which implies that
$\mathbb{P}\{ \tau=\infty\}>0$.
We note that $(\tilde Z_1, \tilde Z_2)$ is a
two-dimensional SRBM with data $(\mathbb{R}_+^2, \hat\theta, \hat
\Gamma,
\hat R)$, where $\hat
\theta=(\theta_1, \theta_2)'$ and $\tilde\Gamma$ is the $2\times2$
principal submatrix of $\Gamma$ obtained by deleting the 3rd row and
the 3rd
column of $\Gamma$. Lemma~\ref{lem:2dStability} in
Appendix \hyperref[sec:AppendixD]{D} will show that when $\hat R$ is a $\mathcaligr{P}$-matrix and the condition (\ref{eq:LCPv2}) is satisfied, the
two-dimensional SRBM $\tilde Z$ is ``rate stable'' in the sense that
%
\begin{equation}
\label{eq:Z12SLLN}
\lim_{t\to\infty} \frac{1}{t}\tilde Z_i(t) =0\qquad  \mbox{almost surely, }
  i=1, 2.
\end{equation}

Solving for $\tilde Y_1$ and $\tilde Y_2$ in the first two components
of (\ref{11}) and plugging them into the third component yields
\begin{eqnarray}
\label{eq:715}
\hspace*{19pt}\tilde Z_3(t)&=&N+ \theta_3 t + \tilde B_3(t)\hspace*{-19pt}\nonumber\\[-8pt]\\[-8pt]
\hspace*{19pt}&&{} + (b, b')\tilde R^{-1}
\biggl[
\pmatrix{ \tilde Z_1(t) \cr
\tilde Z_2(t)
}
-
\pmatrix{
\theta_1 \cr \theta_2
}
t
-
\pmatrix{
\tilde B_1(t)\cr
\tilde B_2(t)
}
\biggr] \tilde Y_3(t), \qquad   t\ge0.\hspace*{-19pt}\nonumber
\end{eqnarray}
%
Equations~(\ref{eq:LCPv3}), (\ref{eq:Z12SLLN}) and (\ref{eq:715}),
together the SLLN for Brownian motion,
imply that
\[
\liminf_{t\to\infty} \frac{1}{t} \tilde Z_3(t) \ge v_3\qquad  \mbox{almost surely}.
\]
Because $v_3>0$, one can argue as in (\ref{eq:78a}) and (\ref{eq:78b})
that for $N$ large enough,
$ \mathbb{P}\bigl\{ \tilde Z_3(t)>1 \mbox{ for all } t\ge0 \bigr\}>0$.
This proves
(\ref{eq:tauhat}).
\end{pf}

Before stating and proving Lemma \ref{lem:caseV} for a LCP solution
in Category \ref{cat5}, we state the following lemma, which is needed
in the proof of Lemma \ref{lem:caseV} and will be proved at the end
of this section.

\begin{lemma}\label{lem:Einfinite}
Let $B=(B_1, B_2, B_3)$ be a three-dimensional Brownian motion with
zero drift and covariance matrix
$\Gamma$, starting from $0$. 
Set
%
\begin{eqnarray}
 Z_1(t) &=& - t + B_1(t) + Y_1(t), \label{eq:Lemma10Z1} \\
 Z_2(t)&=& 2 + B_2(t) - B_1(t)+ Z_1(t), \label{eq:Lemma10Z2} \\
 Z_3(t) &=& 4N + 3\mu t + B_3(t) + a B_1(t)-a
Z_1(t) \label{eq:Lemma10Z3}
\end{eqnarray}
for $t\ge0$, and given constants $a$, $\mu>0$ and $N\ge1$,
where
\[
Y_1(t)= \sup_{0\le s\le t}\bigl(- s +
B_1(s)\bigr)^-.
\]
Then for sufficiently large $N$, one
has
$\mathbb{E}(
\sigma)=\infty$, where $\sigma=\inf\{t\ge0\dvtx  Z_2(t)=1$ or $Z_3(t)=1\}$.
\end{lemma}


\begin{lemma} \label{lem:caseV} If there is a LCP solution
$(u, v)$ in Category \textup{\ref{cat5}}, then $Z$ is not positive recurrent.
\end{lemma}

\begin{pf}
Without loss of generality, we assume that
%
\begin{equation}
\label{eq:uvII}
\qquad \quad  u_1>0, \qquad   u_2=0,\qquad    u_3=0,\qquad    v_1=0, \qquad   v_2= 0, \qquad   v_3>0.
\end{equation}
Then a minor variation of the first paragraph in the proof of Lemma
\ref{lem:theta} establishes that both $\theta_1$ and $\theta_2$ are
negative, so $\theta_1=\theta_2=-1$ with our scaling convention.
Assuming $R$ is as in (\ref{eq:Rspecial}), it
follows as in (\ref{eq:specialLCP}) that $u_1=1$, $a=1$ and
$v_3=b+\theta_3>0$.

Let $Z(0)=(0, 2, N)'$ for some constant $N>1$ and let
\[
\tau_i =\inf\{t\ge0\dvtx  Z_i(t)=1\}, \qquad   i=2, 3,
\]
with $\tau=\min(\tau_2, \tau_3)$.
We will show that $\mathbb{E}_x(\tau)=\infty$ for
sufficiently large $N$, which implies that $Z$ is not
positive recurrent.

The~SRBM $Z$
satisfies equations (\ref{eq:Z1})--(\ref{eq:Z3}). Since $Z_2(t)>0$ and
$Z_3(t)>0$ for $t<\tau$, one has $Y_2(t)=Y_3(t)=0$ for $t<\tau$.
Because $a=1$, (\ref{eq:Z1})--(\ref{eq:Z3})
reduce to
%
\begin{eqnarray}
 Z_1(t) &=& - t + B_1(t) + Y_1(t), \label{eq:IIIZ1} \\
 Z_2(t) &=& 2 - t + B_2(t) + Y_1(t),\label{eq:IIIZ2}\\
 Z_3(t) &=& N + \theta_3 t + B_3(t) + b Y_1(t) \label{eq:IIIZ3}
\end{eqnarray}
on $t< \tau$.
By (\ref{eq:IIIZ1}), one has $Y_1(t)=Z_1(t)+t - B_1(t)$ for
$t<\tau$.
Substituting $Y_1(t)$ into (\ref{eq:IIIZ2}) and (\ref{eq:IIIZ3}), one has
\begin{eqnarray*}
 Z_2(t)&=&2 + B_2(t) - B_1(t)+ Z_1(t), \\
 Z_3(t)& =& N + v_3 t + B_3(t) -b B_1(t)+ b Z_1(t)
\end{eqnarray*}
on $t<\tau$.

For each $t\ge0$, let $\hat Y_1(t)= \sup_{0\le s\le t}(- s +
B_1(s))^-$, and set
\begin{eqnarray*}
 \hat Z_1(t) &=& - t + B_1(t) + \hat Y_1(t),\\
 \hat Z_2(t)&=&2 + B_2(t) - B_1(t)+ \hat Z_1(t), \\
\hat Z_3(t) & =& N + v_3t + B_3(t) -b B_1(t)+ b
\hat Z_1(t)
\end{eqnarray*}
for $t\ge0$.
Let $\hat\tau=\inf\{t\ge0\dvtx  \hat Z_2(t)=1$ or $\hat Z_3(t)=1\}$; clearly, $\tau=\hat\tau$
on every sample path. It follows from Lemma \ref{lem:Einfinite} that
$\mathbb{E}(\hat
\tau)=\infty$ for sufficiently large $N$.
Therefore, $\mathbb{E}_x(\tau)=\infty$, and so
$Z$ is not positive recurrent.
\end{pf}

\begin{pf*}{Proof of Lemma~\ref{lem:Einfinite}}
We first prove the case when $a>0$. When $a \le0$, the proof
is actually significantly simpler; an outline for this
case will be presented at the end of this proof.

Let $X_2(t)=1+B_2(t)-B_1(t)$ and $X_3(t)=B_3(t)+a B_2(t)$ for $t\ge
0$. Then
$X_2$ is a Brownian motion starting from $1$, $X_3$ is a Brownian
motion starting from $0$, and
(\ref{eq:Lemma10Z2})--(\ref{eq:Lemma10Z3}) become
%
\begin{eqnarray}
 Z_2(t)&=&1+ X_2(t)+ Z_1(t)\ge1+ X_2(t), \label{eq:Z2special}\\
\label{eq:Z3special} Z_3(t) &=& (N+a)+ \bigl(N + \mu t + X_3(t)\bigr)\nonumber\\[-8pt]\\[-8pt]
 &&{} + \bigl(2N+2\mu t
-aX_2(t)- a Z_1(t)\bigr)\nonumber
\end{eqnarray}
for $t\ge0$. Define
\begin{eqnarray*}
 \tau_1&=& \inf\{t\ge0\dvtx  X_2(t)\le0\}, \\
 \tau_2 &=& \inf\{t\ge0\dvtx a X_2(t)\ge2N+2\mu t -a Z_1(t) \mbox{ or }
a X_2(t)\ge N + \mu t\},
\\
 \tau_3 &=& \inf\{t\ge0\dvtx  X_3(t) \le-( N+\mu t)\},\\
 \tau&=& \tau_1 \wedge\tau_2\wedge\tau_3.
\end{eqnarray*}
Assuming $N+a>1$, it follows from the definition of $\tau$ and
(\ref{eq:Z2special})--(\ref{eq:Z3special}) that
for each $t<\tau$, $Z_2(t)>1$ and $Z_3(t)>1$, and so $\tau\le
\sigma$. To prove the lemma, it therefore suffices to show that
%
\begin{equation}
\label{eq:tauinfty}
\mathbb{E}(\tau)=\infty.
\end{equation}

Since $X_2$ is a driftless Brownian motion,
$\mathbb{E}(\tau_1)=\infty$. When $N$ is large, it is intuitively clear
that $\tau_1>\tau_2\wedge\tau_3$ with only negligible probability, which
leads to $\mathbb{E}(\tau)=\infty$. To make the argument rigorous,
first note that,
because $Z_1$ is adapted to $B_1$, each $\tau_i$ is a stopping time
with respect to the filtration generated by the Brownian motion $B$,
and hence $\tau$ is a stopping time as well.
Because $X_2(0)=1$ and $X_2$ is a
martingale with respect to the filtration generated by $B$,
by the optional sampling theorem,
\[
\mathbb{E}\bigl(X_2(\tau\wedge t)\bigr)=1
\]
for each $t\ge0$.
We will show that, for sufficiently large $N$,
%
\begin{equation}
\label{eq:half}
\mathbb{E}\bigl(X_2(\tau)1_{\{\tau<\infty\}}\bigr)\le\tfrac{1}{2}.
\end{equation}
It follows that
\begin{eqnarray*}
1&=& \mathbb{E}\bigl(X_2(\tau\wedge t)\bigr)\\
 &=&
\mathbb{E}\bigl(X_2(\tau)1_{\{\tau<t\}}\bigr) + \mathbb{E}\bigl(X_2(t)1_{\{\tau
\ge t\}}\bigr) \\
&\le& \mathbb{E}\bigl(X_2(\tau)1_{\{\tau<\infty\}}\bigr) +
\mathbb{E}\bigl
(\bigl((N+\mu
t)/a\bigr) 1_{\{\tau\ge t\}}\bigr) \\
&=& \tfrac{1}{2} + \bigl((N+\mu t)/a\bigr)\mathbb{P}\{\tau\ge t\},
\end{eqnarray*}
where we have used $X_2(t)\le(N + \mu t)/a$ on $t\le
\tau_2$ for the inequality.
Consequently,
%
\begin{equation}
\label{eq:tauEstimate}
\mathbb{P}\{\tau\ge t\}\ge\frac{a}{2(N+\mu t)}
\end{equation}
for each $t\ge0$,
from which $\mathbb{E}(\tau)=\infty$ follows.

It remains to prove (\ref{eq:half}).
Because $X_2(\tau_1)=0$ when $\tau_1$ is finite,
%
\begin{eqnarray} \label{eq:718}
\mathbb{E}\bigl(X_2(\tau)1_{\{\tau<\infty\}}\bigr)
&=& \mathbb{E}\bigl(X_2(\tau_2\wedge
\tau_3)1_{\{\tau_2\wedge\tau_3< \tau_1\}}\bigr) \nonumber\\
&\le& \mathbb{E}\bigl(\bigl(\bigl(N+ \mu(\tau_2\wedge\tau_3)\bigr)/a\bigr)1_{\{\tau
_2\wedge
\tau_3<
\infty\}}\bigr) \\
&\le& \sum_{n=0}^\infty\bigl(\bigl(N+\mu(n+1)\bigr)/a\bigr) \mathbb{P}\{n <
\tau_2\wedge\tau_3\le
n+1\}.\nonumber
\end{eqnarray}

To bound the probability $\mathbb{P}\{n < \tau_2\wedge\tau_3\le
n+1\}$ for each $n\in\mathbb{Z}_+$, we use
\[
\{n < \tau_2\wedge\tau_3\le
n+1\}\subset
\{n < \tau_2\le
n+1\} \cup\{n < \tau_3\le
n+1\}.
\]
For $\{n < \tau_3\le n+1\}$, $X_3(t)\le
-(N+\mu n)$ first occurs on $(n, n+1]$. By the
strong Markov property for Brownian motion and the
reflection principle, the
probability of the latter event is at most $2\mathbb{P}\{X_3(n+1)<-( N
+ \mu n) \}$. For $\{n < \tau_2\le n+1\}$, either
$aX_2(t) \ge N+\mu n$ first occurs on $(n, n+1]$,
or $aZ_1(t) \ge N+\mu n$ occurs on $(n,
n+1]$.
One can also apply the strong Markov property and the reflection
principle to the first event. One therefore obtains
\begin{eqnarray}\label{eq:727}
\mathbb{P}\{n < \tau_2\wedge\tau_3\le
n+1\}
&\le& \mathbb{P}\{n < \tau_2\le n+1\} + \mathbb{P}\{n < \tau_3\le
n+1\}
\nonumber\\
&\le& 2 \mathbb{P}\{X_3(n+1)<-(N + \mu n)\}\nonumber\\[-8pt]\\[-8pt]
&&{} + 2\mathbb{P}\{
aX_2(n+1)> N + \mu
n\}\nonumber \\
&& {}+ \mathbb{P}\Bigl\{ \sup_{n< s \le n+1}aZ_1(s)> N + \mu n
\Bigr\}.\nonumber
\end{eqnarray}

We proceed to bound each of these three terms and plug these bounds
into the last term in (\ref{eq:718}).
Note that
%
\begin{eqnarray}\label{eq:731a}
\mathbb{P}\{ X_3(n+1)< -(N+ \mu n)\}
& = & \mathbb{P}\biggl\{N(0, 1)>
\frac{N+\mu n}{\gamma\sqrt{n+1}}\biggr\} \nonumber\\
&\le& \frac{1}{\sqrt{2\pi}} \frac{\gamma\sqrt{n+1}}{N+\mu n}
\exp\biggl(-\frac{1}{2}\frac{(N+\mu n)^2}{\gamma^2(n+1)}\biggr)
\\
&\le& \frac{\gamma}{\sqrt{2\pi}\mu}
\exp\biggl(-\frac{\mu}{2\gamma^2}(N+\mu n)\biggr) \nonumber
\end{eqnarray}
for $N\ge\mu$,
where $N(0, 1)$ denotes the standard normal random variable and
$\gamma^2$ is the variance of the Brownian motion $X_3$.
The~first inequality is a standard estimate and is obtained by
integrating by parts. 
Assume $N=N' \mu$, with $N'\in\mathbb{Z}_+$. Then,
\begin{eqnarray}\label{eq:728}
&& \sum_{n=0}^\infty\bigl(N+\mu(n+1)\bigr) 2 \mathbb{P}\{ X_3(n+1)<-(N+
\mu n)\}
\nonumber\\[-8pt]\\[-8pt]
&&\qquad  \le\sum_{n=N'}^\infty\frac{2\gamma}{\sqrt{2\pi}}(n+1)
\exp\biggl(-\frac{\mu^2}{2\gamma^2} (n+1)\biggr),\nonumber
\end{eqnarray}
which is less than $1/6$ for sufficiently large $N$. For the same reason,
%
\begin{eqnarray}
\sum_{n=0}^\infty\bigl(N+\mu(n+1)\bigr) 2 \mathbb{P}\{a X_2(n+1)> N+ \mu n\}
\le\frac{1}{6}\label{eq:719}
\end{eqnarray}
for sufficiently large $N$.
To bound the probability $\mathbb{P}\bigl\{ \sup_{n< s \le n+1}aZ_1(s)>
N + \mu n\bigr\}$, we apply Lemma \ref{lem:rbmEstimate} in
Appendix~\hyperref[sec:AppendixD]{D}. The~lemma states that for appropriate
constants $c_1,
c_2>0$,
\begin{eqnarray*}
\mathbb{P}\Bigl\{\sup_{n\le s \le n+1}Z_1(t)>x\Bigr\} \le c_1 \exp(-c_2
x)\qquad
\mbox{for } x\ge0.
\end{eqnarray*}
For $N=N'\mu$,
\begin{eqnarray}\label{eq:729}
&& \sum_{n=0}^\infty\bigl(N+\mu(n+1)\bigr)\mathbb{P}\Bigl\{ \sup
_{n< s
\le n+1}aZ_1(s)>
N + \mu n\Bigr\} \nonumber\\[-8pt]\\[-8pt]
&&\qquad  \le\sum_{n=N'}^\infty(n+1)\mu c_1 \exp(- c_2\mu n/a),
\nonumber
\end{eqnarray}
which is also less than $1/6$ for large $N$. The~bounds obtained for
(\ref{eq:728})--(\ref{eq:729}) together show that the last term in
(\ref{eq:718}) is at most $1/2$. This implies (\ref{eq:half}), and
hence the lemma for $a>0$.

When $a\le0$, the proof is analogous to the case $a>0$, with the
following simplifications. Assume $N+a\ge1$. The~equality
(\ref{eq:Z3special}) can be replaced by
\[
Z_3(t) \ge(N+a)+ \bigl(N + \mu t + X_3(t)\bigr) + \bigl(2N+2\mu t
-aX_2(t)\bigr)
\]
for $t\ge0$, because $a\le0$. The~definition of $\tau_2$ can be
replaced by the simpler
\[
\tau_2 = \inf\{t\ge0\dvtx  X_2(t)\ge N+\mu t \}.
\]
One can again check that for each $t<\tau$, $Z_2(t)>1$ and
$Z_3(t)>1$, and so $\tau\le\sigma$. To prove the lemma for the case
$a\le0$, it remains to show
(\ref{eq:tauinfty}). For this, we follow the same procedure as in the
case $a>0$. First, (\ref{eq:half}) still implies
(\ref{eq:tauEstimate}) with $a$ in the right side of
(\ref{eq:tauEstimate}) replaced by $1$. From (\ref{eq:tauEstimate}),
(\ref{eq:tauinfty}) follows. To demonstrate (\ref{eq:half}), we employ
(\ref{eq:718}), with $a$ there replaced by $1$.
To bound the probability $\mathbb{P}\{n < \tau_2\wedge\tau_3\le
n+1\}$ for each $n\in\mathbb{Z}_+$, we replace (\ref{eq:727}) with the
simpler
\begin{eqnarray*}
&&\mathbb{P}\{n < \tau_2\wedge\tau_3\le
n+1\}\\
&&\qquad \le \mathbb{P}\{n < \tau_2\le n+1\} + \mathbb{P}\{n < \tau_3\le
n+1\}
\nonumber\\
&&\qquad \le 2 \mathbb{P}\{X_3(n+1)<-(N + \mu n)\} + 2\mathbb{P}\{X_2(n+1)>
N + \mu
n\} . \label{eq:727a}
\end{eqnarray*}
It then follows from bounds (\ref{eq:731a})--(\ref{eq:719}) that
(\ref{eq:half}) holds for sufficiently large~$N$. This implies the
lemma for $a\le0$.
\end{pf*}

\begin{appendix}

\section{Semimartingale reflecting Brownian motions}
\label{sec:srbm}

In this section, we present the standard definition of a
semimartingale reflecting Brownian motion (SRBM) in the
$d$-dimensional orthant $S=\mathbb{R}^d_+$, where $d$ is a positive
integer. We also review the standard definition of positive
recurrence for an SRBM, connecting it with the alternative
definition used in Section \ref{sec1}.

Recall from Section \ref{sec1} that $\theta$ is a constant vector in
$\mathbb{R}^d$, $\Gamma$ is a $d\times d$ symmetric and strictly positive
definite matrix, and $R$ is a $d\times d$ matrix. We shall define an
SRBM associated with the data $(S, \theta, \Gamma, R)$. For
this, a
triple $(\Omega, \mathcaligr{F}, \{\mathcaligr{F}_t\} )$ will be called a
\textit{filtered space} if $ \Omega$ is a set, $ \mathcaligr{F}$ is a
$\sigma$-field of subsets of $\Omega$, and $\{\mathcaligr{F}_t\} \equiv
\{
\mathcaligr{F}_t, t \geq0\}$ is an increasing family of sub-$\sigma
$-fields of
$\mathcaligr{F}$, that is, a filtration.

\begin{definition}[(Semimartingale reflecting Brownian motion)] \label{def:srbm}
A SRBM associated with $(S, \theta, \Gamma, R)$ is a continuous
$\{\mathcaligr{F}_t\}$-adapted $d$-dimensional process $Z=\{Z(t), t\ge
0\}$, together with a family of probability measures $\{\mathbb{P}_x, x
\in S\}$,
defined on some filtered space $(\Omega, \mathcaligr{F}, \{\mathcaligr
{F}_t\})$ such that,
for each $x\in S$, under~$\mathbb{P}_x$, \eqref{11} and \eqref{14} hold,
where, writing $W(t)=X(t)-\theta t$ for $t\ge0$, $W$ is a
\mbox{$d$-dimensional} Brownian motion with
covariance matrix $\Gamma$, an $\{\mathcaligr{F}_t\}$-martingale such
that \mbox{$W(0)=x$}
$\mathbb{P}_x$-a.s., and $Y$ is an
$\{\mathcaligr{F}_t\}$-adapted $d$-dimensional process such that
$\mathbb{P}_x$-a.s.  \eqref{12} and \eqref{13} hold. Here \eqref
{12} is
interpreted
to hold for each component of $Y$, and \eqref{13} is defined to be
%
\begin{equation}
\label{eq:1}
\int_0^t 1_{\{Z_i(s)\neq0\}}\,  d Y_i(s)=0\qquad   \mbox{for all } t\ge0.
\end{equation}
\end{definition}

Definition~\ref{def:srbm} gives the so-called weak formulation of a
SRBM. It is a standard definition adopted in the literature; see,
for example, \citet{dupwil94} and \citet{wil95}. Note that
condition \eqref{eq:1} is equivalent to the condition that, for
each $t>0$, $Z_j(t)>0$ implies $Y_j(t-\delta)=Y_j(t+\delta)$ for
some $\delta>0$. \citet{reiwil88} showed that a necessary condition
for a $(S, \theta, \Gamma, R)$-SRBM to exist is that the
reflection matrix $R$ is completely $\mathcaligr{S}$ (this term was
defined in Section \ref{sec1}). \citet{taywil93} showed that when $R$ is
completely $\mathcaligr{S}$, a~$(S, \theta, \Gamma, R)$-SRBM $Z$
exists and $Z$ is unique in law under $\mathbb{P}_x$ for each $x\in
S$. Furthermore, $Z$, together with the family of
probability measures $\{\mathbb{P}_x, x\in\mathbb{R}^d_+\}$, is a Feller
continuous strong Markov process.

Let $(\theta, \Gamma, R)$ be fixed with $\Gamma$ being a positive
definite matrix and $R$ being a completely $\mathcaligr{S}$ matrix.
\citeauthor{dupwil94} [(\citeyear{dupwil94}), Definition 2.5] and
\citeauthor{wil95} [(\citeyear{wil95}), Definition 3.1] gave
the following definition of
positive recurrence.

\begin{definition}
\label{def:positive-recurrence}
An SRBM $Z$ is said to be \textit{positive recurrent} if, for each
closed set
$A$ in $S$ having positive Lebesgue measure, we have
$\mathbb{E}_x(\tau_A)<\infty$ for all $x\in S$, where
$\tau_A=\inf\{t\ge0\dvtx  Z(t)\in A\}$.
\end{definition}

Because each open neighborhood of the origin contains a closed
ball that has positive volume, Definition
\ref{def:positive-recurrence} appears to be stronger (i.e.,
more restrictive) than the definition adopted in Section \ref{sec1}, but
one can show that these two notions of positive recurrence are
equivalent for a SRBM. Indeed, the last paragraph on page 698 in
\citet{dupwil94} provides a sketch of that proof.

\section{Convenient normalizations of problem data}
\label{sec:normalization}

Let $R$ be a $d\times d$ completely $\mathcaligr{S}$ matrix and $(X, Y, Z)$ a
triple of continuous, \mbox{$d$-dimensional} stochastic processes defined on
a common probability space. The~diagonal elements of $R$ are
necessarily positive. Let $\tilde D =\operatorname{diag}(R)$ and $\tilde
R = R \tilde D^{-1}$ (thus $\tilde R$ is a $d\times d$
completely-$\mathcaligr{S}$
matrix that has ones on the diagonal), and define $\tilde Y(t) =
\tilde D Y (t)$ for $t \ge0$. If ($X$, $Y$, $Z$) satisfy
(\ref{11})--(\ref{14})
with reflection matrix $R$, then ($X$, $\tilde Y$, $Z$) satisfy
(\ref{11})--(\ref{14})
with reflection matrix $\tilde R$, and vice versa. Thus the
distribution of $Z$ is not changed if one substitutes $\tilde R$ for
$R$, and that substitution assures the standardized problem format
(\ref{eq:new1.5}).

Now let $R$ and ($X$, $Y$, $Z)$ be as in the previous paragraph, and
further suppose that $X$ is a Brownian motion with drift vector
$\theta$ and nonsingular covariance matrix~$\Gamma$. Define a
$d\times d$ diagonal matrix $D$ by setting $D_{ii} = 1$ if $\theta_i =
0$ and\vspace*{1pt} $D_{ii} = \vert\theta_i\vert^{-1}$ otherwise. Setting $\hat Z
= DZ$, $\hat X = DX$ and $\hat R = DR$, one sees that if ($X$, $Y$,
$Z$) satisfy (\ref{11})--(\ref{14}) with reflection matrix $R$, then
($\hat X$, $Y$, $\hat Z$) satisfy (\ref{11})--(\ref{14}) with reflection
matrix $\hat R$, and vice versa. Of course, $\hat X$ is a Brownian
motion whose drift vector $\hat\theta= D \theta$ satisfies
(\ref{eq:51}); the covariance matrix of $\hat X$ is $\hat\Gamma= D
\Gamma D$. Thus our linear change of variable gives a transformed
problem in which (\ref{eq:51}) is satisfied. To achieve a problem
format where \textit{both} (\ref{eq:new1.5}) and (\ref{eq:51}) are
satisfied, one can first make the linear change of variable described
in this paragraph, and then make the substitution described in the
previous paragraph.

\section{Proof that (\protect\ref{15}) is necessary for stability of $Z$}
\label{sec:proofNecessa}

We consider a $d$-dimensional SRBM $Z$ with associated data ($S,
\theta, \Gamma, R$), defined as in Appendix \hyperref[sec:srbm]{A}, assuming throughout
that $R$ is completely $\mathcaligr{S}$. Let us also assume until further
notice that $R$ is nonsingular. Because $R$ is an $\mathcaligr{S}$-matrix,
there exist $d$-vectors $w$, $v > 0$ such that $Rw = v$. That is,
%
\begin{equation}\label{c1}
R^{-1}v > 0\qquad  \mbox{where } v >
0.
\end{equation}

Now suppose it is \textit{not} true that
$R^{-1}\theta<0$. That is,
defining $\gamma= R^{-1}\theta$, suppose that $\gamma_{i} \ge$ 0
for some $i\in\{1,\ldots, d\}$. For future reference let $u$ be
the $i${th} row of~$R^{-1}$. Thus (\ref{c1}) implies
%
\begin{equation}\label{c2}
u\cdot v > 0.
\end{equation}

Our goal is to show that $Z$ cannot be positive recurrent. Toward
this end, it will be helpful to represent the underlying Brownian
motion $X$ in (\ref{11}) as $X(t)=W(t)+\theta t$, where $W$ is a
\mbox{$d$-dimensional} Brownian motion with zero drift and covariance
matrix $\Gamma$. Premultiplying both sides of (\ref{11}) by $R^{-1}$
then gives $R^{-1}Z(t)=R^{-1}W(t)+\gamma t+Y(t)$. The~$i${th}
component of the vector equation is
%
\begin{equation}\label{c3}
\xi(t)\equiv u\cdot Z(t)=u\cdot W(t)+\gamma
_{i}t+Y_{i}(t),\qquad    t\ge0.
\end{equation}
Let $A = \{ z\in S \dvtx  \vert z\vert\le1\}$
and $B = \{ u\cdot z \dvtx  z\in A\}$. Then $B\subset\mathbb{R}$ is a
compact interval containing the origin, and from (\ref{c2}) we
know that $B$ contains positive values as well, because $A$
contains $\alpha v$ for sufficiently small constants $\alpha
>$ 0. Thus $B$ has the form
%
\begin{equation}\label{c5} B = [ a , b ]\qquad  \mbox{where } a \le0
\mbox{ and } b >
0.
\end{equation}

As the initial state $x=Z (0) = W (0)$, we take
%
\begin{eqnarray}\label{c6}
x&=&\beta v \qquad  \mbox{where $v$ is chosen as
in (\ref{c1})\quad  and }\nonumber\\[-8pt]\\[-8pt]
 \beta
&>& \max\bigl( \vert v\vert^{-1 }, ( u\cdot v)^{ -1}b\bigr) .\nonumber
\end{eqnarray}
From (\ref{c1}), (\ref{c2}), (\ref{c5}) and (\ref{c6}), we have
that
%
\begin{equation}\label{c7}
x\in S ,\qquad   \vert x\vert> 1\quad
\mbox{and}\quad   u\cdot x > b.
\end{equation}
Thus, defining $\tau_{A} = \inf\{t \ge0\dvtx  Z(t)\in A\}$ and
$\sigma= \inf\{t \ge0\dvtx  \xi(t)\in B\}$, it follows from
the definitions of $A$, $B$ and $\xi$, plus (\ref{c5}) and
(\ref{c7}), that
%
\begin{equation}\label{c8}
\tau_{A} \ge\sigma,\qquad
\mbox{$\mathbb{P}_{x}$-a.s.}
\end{equation}
From (\ref{c3}), we see that $\xi$ is bounded below
by a one-dimensional Brownian motion with nonnegative drift, and
$\xi(0) > b$ $\mathbb{P}_{x}$-a.s. Thus $\mathbb{E}_{x}(\sigma
)=\infty$, implying that $\mathbb{E}_{x}(\tau_{A})=\infty$ as well
by (\ref{c8}). This establishes that $Z$ is not positive
recurrent when $R$ is nonsingular.

We still need to show that $Z$ cannot be positive
recurrent when $R$ is singular. In this case, there exists a
nontrivial vector $u\in\mathbb{R}^d $ such that $u' R$ = 0, and
we can assume that $u\cdot\theta\ge0$ as well (because $-u$ can
be exchanged for $u$ if necessary). Premultiplying both sides of
(\ref{11}) by $u'$ gives the following analog of
(\ref{c3}):
\[
u\cdot Z(t)=u\cdot W(t) + ( u\cdot\theta)t .
\]
Because $R$ is an
$\mathcaligr{S}$-matrix, for this given $u$ there exist $w, v\in S$ such that
%
\begin{equation}\label{c9}
u + Rw = v .
\end{equation}
After premultiplying both sides of (\ref{c9}) by $u^\prime$, one
obtains
\[
u\cdot v=\vert u\vert^{2} >0.
\]
We choose the initial state $x=Z(0) = W(0)$
exactly as in (\ref{c6}), and define the set $A$ as before. The~proof that $\mathbb{E}_{x}(\tau_{A})=\infty$, and hence that $Z$ is
not positive recurrent, then proceeds exactly as in the case
treated above, except that now the process $\xi=u\cdot Z$ is
itself a Brownian motion with nonnegative drift, whereas in the
case treated earlier $\xi$ was bounded below by such a Brownian
motion.

\section{Two lemmas}
\label{sec:AppendixD}

In this appendix, we demonstrate two lemmas that are used in
Section~\ref{sec:lemmas}.
Lemma~\ref{lem:rbmEstimate} is employed in the proof of Lemma
\ref{lem:Einfinite}.

\begin{lemma} \label{lem:rbmEstimate} Let $X$ be a one-dimensional
Brownian motion with drift
$\theta<0$, variance $\sigma^2$, starting from $0$. Let $Z$ be the
corresponding one-dimensional SRBM,
\[
Z(t) = X(t) -\min_{0\le s \le t} X(s)\qquad   \mbox{for } t\ge0.
\]
There exist constants $c_1>0$ and $c_2>0$ such that
%
\begin{equation}
\label{eq:estimates}
\mathbb{P}\Bigl\{\sup_{n-1\le s \le n} Z(s)> x\Bigr\} \le c_1 \exp
(-c_2 x)
\end{equation}
for all $n\in\mathbb{Z}_+$ and $x>0$.
\end{lemma}

\begin{pf}
One could employ the Lipschitz continuity property of the
one-dimensional Skorohod map and the estimate (4.9) of
\citet{atar-2001-29} to prove the lemma. Here, we provide a direct
proof. Since $X$ has negative drift, it is well known [see, e.g.,
Section 1.9 of \citet{har85}] that, for each $t>0$ and
$x>0$,
%
\begin{eqnarray}\label{eq:rbmexp}
\mathbb{P}\{Z(t)>x\} & = & \mathbb{P}\Bigl\{\max_{0\le s \le
t}\bigl(X(t)-X(t-s)\bigr)>x\Bigr\}\nonumber\\
&=& \mathbb{P}\Bigl\{\max_{0\le s \le t} X(s)>x\Bigr\} \\
&\le&\mathbb{P}\Bigl\{\sup_{0\le s < \infty} X(s)>x\Bigr\}= \exp
(-
2|\theta|x/\sigma^2).\nonumber
\end{eqnarray}
Therefore, for $n\in\mathbb{Z}_+$,
\begin{eqnarray*}
&&\mathbb{P}\Bigl\{\max_{n-1< s \le n}Z(s)>x\Bigr\}\\
 &&\qquad \le
\mathbb{P}\{Z(n)>x/2\}+ \mathbb{P}\Bigl\{ Z(n)\le x/2,
\max_{n-1< s \le n}Z(s)>x\Bigr\}\\
&&\qquad \le \exp(-|\theta|x/\sigma^2) +
\mathbb{P}\{ n-1< \tau< n, X(n)-X(\tau)< -x/2\},
\end{eqnarray*}
where $\tau=\inf\{s\ge n-1\dvtx  Z(s)> x\}$.

Note that $\tau$ is a stopping
time with respect to the filtration generated by $Z$, which is
the filtration generated by the Brownian motion $B$, where
$B(t)=X(t)-\theta t$ for $t\ge0$.
By the strong Markov property of $B$,
\begin{eqnarray*}
&&\mathbb{P}\{ n-1< \tau< n, X(n)- X(\tau)< -x/2\}\\
 &&\qquad \le
\mathbb{P}\{ n-1< \tau< n, B(n)- B(\tau)< -x/2+|\theta|\} \\
&&\qquad \le \Phi\bigl((-x/2+|\theta|)/\sigma\bigr)
\end{eqnarray*}
for $x> 2|\theta|$, where $\Phi$ is the standard normal distribution function.
Thus,
\begin{eqnarray*}
\mathbb{P}\Bigl\{\max_{n-1< s \le n}Z(s)>x\Bigr\}
\le\exp(-|\theta|x/\sigma^2) +
\Phi\bigl((-x/2+1)/\sigma\bigr)
\end{eqnarray*}
for $x>2|\theta|$, from which (\ref{eq:estimates}) follows.
\end{pf}

The~following lemma is used in the proof of
Lemma~\ref{lem:caseII}.
Let
\[
\label{eq:D22}
R=
\pmatrix{
1 & b \cr
a & 1
}
\quad  \mbox{and}\quad
\theta=
\pmatrix{
\theta_1 \cr
\theta_2
}
,
\]
and assume that $\det(R)>0$. Then $R$ is a $\mathcaligr{P}$-matrix and hence
is completely $\mathcaligr{S}$. For any given $2\times2$ positive
definite matrix $\Gamma$, it therefore follows from \citet{taywil93}
that, starting from any fixed state $x\in\mathbb{R}_+^2$, there is a
two-dimensional SRBM $Z$ with data $(\theta,\Gamma, R)$, as defined
in Definition \ref{def:srbm}, that is well defined and is unique in
law.

\begin{lemma}\label{lem:2dStability}
Suppose that $R^{-1}\theta\le0$. Then each fluid path $(y, z)$
starting from
$z(0)=0$ remains at $0$; that is, $z(t)=0$ for $t\ge0$. Consequently,
%
\begin{equation}
\label{eq:rateStablity}
\lim_{t\to\infty} \frac{Z(t)}{t} = 0 \qquad  \mbox{almost surely}.
\end{equation}
%
\end{lemma}

\begin{pf}
It follows from
\[
R^{-1}=\frac{1}{1-ab}
\pmatrix{
1 & -b \cr
-a & 1
}
,
\]
$\det(R)=1-ab>0$ and $R^{-1}\theta\le0$ that
\[
\label{eq:D2theta}
\theta_1 - b \theta_2 \le0,\qquad   \theta_2 - a \theta_1 \le0.
\]

Observe that $R^{-1}$ is a $\mathcaligr{P}$-matrix because $R$ is a
$\mathcaligr{P}$-matrix. Therefore, $R^{-1}$ is an $\mathcaligr{S}$-matrix.
Consequently, there exists $(d_1, d_2)>0$ such that $(c_1,
c_2)\equiv(d_1, d_2)R^{-1}>0$. By assumption $R^{-1}\theta\le0$,
and so
one has
\[
\label{eq:ctheta}
c_1 \theta_1 + c_2 \theta_2 = (d_1, d_2)R^{-1}\theta\le0.
\]

Let $(y,z)$ be a fluid path with $z(0)=0$. By the oscillation inequality
[see, e.g., Lemma 4.3 of \citet{daiwil95}], $(y,z)$ is Lipschitz
continuous.
Setting $g(t)=c_1z_1(t)+c_2 z_2(t)$, we will show that $g(t)=0$
for $t\ge0$. From this, it follows that $z(t)=0$ for $t\ge0$, which
is the first claim in the lemma.
It suffices to prove that
\begin{eqnarray*}
\label{eq:ggdot}
\dot g(t) \le0 \mbox{ at each $t$ where } g(t)>0
\mbox{ and $(y, z)$ is differentiable.}
\end{eqnarray*}
We consider several cases, depending on whether $z_1(t)$ and $z_2(t)$
are strictly positive. When $z_1(t)>0$ and $z_2(t)>0$,
(\ref{eq:fl-1})--(\ref{eq:fl-4}) imply that $\dot z_1(t)=\theta_1$ and
$\dot z_2(t)=\theta_2$, and hence $\dot g(t)=c_1 \dot z_1(t)+c_2 \dot
z_2(t)\le0$. When $z_1(t)=0$ and $z_2(t)>0$, $\dot z_1(t)=0$ and
$\dot y_2(t)=0$, from which one has $\dot y_1(t)= - \theta_1$ and
$\dot z_2(t)=\theta_2 -a \theta_1\le0$. Because in this case $\dot
z_1(t)=0$ and
$\dot z_2(t)\le0$, $\dot g(t)\le0$ follows. When $z_2(t)=0$
and $z_1(t)>0$, one can similarly argue that $\dot g(t)\le0$.
This shows that $g(t)=0$ for $t\ge0$, as desired.

We next demonstrate (\ref{eq:rateStablity}).
Let $Z$ be a two-dimensional SRBM with data $(\theta, \Gamma, R)$
having a given
initial point $Z(0)=x\in\mathbb{R}_+^2$. By Definition~\ref
{def:srbm}, $Z$,
together with the associated pair $(X,Y)$, satisfies (\ref{11})--(\ref{14}).
For each $r>0$ and each $t\ge0$, set
\[
\bar{X}^r(t)= \frac{1}{r} X(rt),\qquad
\bar{Y}^r(t)= \frac{1}{r} Y(rt),\qquad
\bar{Z}^r(t)= \frac{1}{r} Z(rt).
\]
By the functional strong law of large numbers (FSLLN) for Brownian motion,
almost surely,
%
\begin{equation}
\label{eq:Xslln}
\lim_{r\to\infty} \sup_{0\le s\le t}|\bar{X}^r(s)- x (s)|=0
\qquad  \mbox{for each } t>0,
\end{equation}
where $x(s)= (\theta_1s, \theta_2 s)'$ for $s\ge0$. Fix a sample
path that satisfies (\ref{eq:Xslln}) and let $\{r_n\}\subset\mathbb
{R}_+$ be
a sequence with $r_n\to\infty$. The~FSLLN (\ref{eq:Xslln}) implies
that $\{\bar{X}^{r_n}\}$ is relatively
compact in $C(\mathbb{R}_+, \mathbb{R})$, the space of continuous
functions on $\mathbb{R}_+$
endowed with the topology of uniform convergence on compact
sets. By the oscillation inequality, $\{(\bar{Y}^{r_n},
\bar{Z}^{r_n})\}$ is also relatively compact.

Let $(y,z)$ be a limit point
of $\{(\bar{Y}^{r_n}, \bar{Z}^{r_n})\}$. It is not difficult to show
that $(y,z)$ is a
fluid path associated with the data $(\theta, R)$ that satisfies
$z(0)=0$. It follows from the first part of the proof that
$z(t)=0$ for $t\ge
0$. This is the unique limit point with $z(0)=0$.
Therefore, almost surely,
%
\[
\lim_{r\to\infty} \frac{1}{r}Z(r)= \lim_{r\to\infty}\bar{Z}^r(1)=0,
\]
which proves (\ref{eq:rateStablity}).
\end{pf}

\end{appendix}

\section*{Acknowledgments}
We thank Ruth Williams, Kavita Ramanan and Richard Cottle for helpful
discussions. In particular, the proof that (\ref{15}) is necessary for
positive recurrence of $Z$ (see Appendix \hyperref[sec:proofNecessa]{C}) is
based on a general line of argument suggested by Ruth Williams.

%

\printaddresses

\end{document}